\documentclass[11pt,leqno]{article}
\usepackage{amsfonts,latexsym,amssymb,amsmath,color}
\usepackage[pdftex]{graphicx}
\def\bn{\mathbf n}
\def\bb{\mathbf b}
\def\bx{\mathbf x}
\def\by{\mathbf y}
\def\bm{\mathbf m}
\def\bk{\mathbf k}
\def\be{\mathbf e}
\def\bP{\mathbf p}
\def\bv{\mathbf v}
\def\bw{\mathbf w}
\def\ba{\mathbf a}
\def\bu{\mathbf u}
\def\bz{\mathbf z}
\def\bq{\mathbf q}
\def\bh{\mathbf h}
\def\bof{\mathbf f}
\def\cL{{\cal L}}
\def\cA{{\cal A}}

\def\cM{{\cal M}}
\def\cS{{\cal S}}
\def\cO{{\cal O}}
\def\cB{{\cal B}}
\def\cH{{\cal H}}

\def\cD{{\cal D}}
\def\cU{{\cal U}}

\def\fA{{\mathfrak A}}
\def\fr{{\mathfrak r}}
\def\fp{{\mathfrak p}}
\def\fa{{\mathfrak a}}
\def\fb{{\mathfrak b}}
\def\fc{{\mathfrak c}}
\def\fd{{\mathfrak d}}
\def\bR{{\Bbb R}}
\def\bC{{\Bbb C}}
\def\bS{{\Bbb S}}
\def\var{\varepsilon}
\newtheorem{proposition}{Proposition}

\newtheorem{definition}[proposition]{Definition}

\numberwithin{equation}{section}

\begin{document}


\title{Analysis of Nematic Liquid Crystals with Disclination Lines}

\author{Patricia Bauman\thanks{Research supported by NSF grants DMS-0456286 and DMS-0604839.}
and
 Daniel Phillips\thanks{Research supported by NSF grants DMS-0456286 and DMS-0604839.}\\
 Department of Mathematics\\
 Purdue University\\
West Lafayette, IN 47906\\
bauman@math.purdue.edu, phillips@math.purdue.edu
\and
 Jinhae Park\thanks{Research supported by NSF grant DMS-0604839.}\\
  Department of Mathematics\\
  Chungnam National University\\
  220 Kung-Dong, Yuseong-Gu\\
  Daejeon 305-763, South Korea,\\
  jhpark2003@gmail.com}
\bibliographystyle{siam}
\date{\today}
\maketitle

\begin{abstract}
We investigate the structure of nematic liquid crystal  thin films
described by the Landau--de Gennes tensor-valued order parameter
with Dirichlet boundary conditions of nonzero degree. We prove that
as the elasticity constant goes to zero a limiting uniaxial texture
forms with disclination lines corresponding to a finite number of
defects, all of degree $1\over 2$ or all of degree $-{1\over 2}$. We
also state a result on the limiting behavior of minimizers of the
Chern-Simons-Higgs model without magnetic field that follows from a
similar proof.
\end{abstract}

\section{Introduction}
We investigate disclination line defects in a nematic liquid crystal
by using a tensor-valued order parameter description based on the
Landau-de Gennes theory. The unknown field $Q$ in this theory is
$\cS$-valued such that $Q=Q(x,y)$, where $\cS$ is the space of
$3\times 3$, real
symmetric, traceless matrices, and $(x,y)$ varies in a bounded
domain $\Omega$ in $\mathbb R^2$.  For simplicity, we assume that
$\Omega$ is a simply connected bounded domain with a $C^3$ boundary
in the plane, representing the reference configuration of a very
thin liquid crystal material.

The Landau-de Gennes model is based on a phenomenological theory in
which stable states of the liquid material correspond to minimizers
(or stable states) of an energy formulated in terms of $Q$ on
$\Omega$.  The matrix $Q(\bx)$ models the second moments of the
orientations of the rod-like liquid crystal molecules near $\bx$.
Its values describe the average orientation and phase of the liquid
crystals near $\bx$, measured through its eigenvectors and
eigenvalues. (See Section 1.1 for more detail on this structure.) As
such $Q$ is a measure of the microscopic anisotropy of their
relative positions.  In this paper, we consider fields $Q\in
W^{1,2}(\Omega;\cS)$ with fixed uniaxial nematic boundary conditions
of the form $Q=Q_0$ on $\partial \Omega$ (in the sense of trace). We
assume throughout the paper that $(Q_0)_{ij} \in
C^3(\partial\Omega)$ for all $1\leq i,j\leq 3$, and
\begin{equation}
Q_0(x,y)=s(\bn_0(x,y)\otimes\bn_0 (x,y)-{1\over 3}\ I)\quad\text{ for }(x,y)\in\partial\Omega
\end{equation}
where {$I$ is the $3\times 3$ identity matrix, $s$ is an arbitrary fixed
nonzero real number, and
$\bn_0$ is a fixed vector field defined on $\partial \Omega$
satisfying $\bn_0=\langle n_1,n_2,0\rangle$, $|\bn_0|=1$, and (1.1) on $\partial
\Omega$.  Note that $Q_0$ is invariant under changes in direction:
$\bn_0(x,y) \to -\bn_0(x,y)$ at any point $(x,y)$ in $\partial \Omega$, which allows
boundary conditions of degree one-half, or integer multiples of
one-half, for $Q_0$.  Nonzero boundary conditions of this type are
observed in thin liquid crystal materials exhibiting defects along
curves, known as "disclination lines," whose cross-sections in
$\Omega$ are isolated points.  (See Figure 1.) We analyze a class of
equilibria for the Landau-de Gennes energy
\begin{equation*}
F_\var(Q)=\int_\Omega [f_e(Q)+\var^{-2} f_b(Q)].
\end{equation*}
where $\var > 0$, defined for all $Q \in W^{1,2}(\Omega,\cS)$.  Here
$f_e$ is the elastic energy density in $\Omega$ given by
\begin{eqnarray*}
f_e(Q)&=&{L_1\over 2}\ Q_{ij,k}\ Q_{ij,k}+{L_2\over 2}\ Q_{ij,j}\ Q_{ik,k}\nonumber\\
&+&{L_3\over 2}\ Q_{ij,k}\ Q_{ik,j},
\end{eqnarray*}
where each term above is summed over all $i,j,k$ from 1 to 3,
$Q_{ij,\alpha}$ denotes $\frac {\partial Q_{ij}} {\partial x_\alpha}$, and
$(x_1,x_2,x_3)=(x,y,z)$.  The above formula is valid in two or
three-dimensional reference domains.  Since here we are considering a
two-dimensional
reference domain $\Omega$, we identify $Q(x,y)$ with $Q(x,y,0)$ above,
so that $Q_{ik,3}=0$ for all $1 \leq i,j \leq 3.$  We assume
throughout the paper that
\begin{equation}
L_1>0\text{ and }L_1+L_2+L_3>0.
\end{equation}
The term $f_b$ is the bulk energy density given by a real-valued
$C^\infty$ function which depends on temperature as well as on $Q$.
We assume that temperature is fixed and $f_b=f_b(Q)$ is a
nonnegative $C^{\infty}$ function defined on $\cS$ such that
$f_b(Q)=0$ if and only if $Q\in\Lambda_s=\{Q\in \cS\colon Q=
s(\bm\otimes\bm-{1\over 3}I)$ for some $\bm\in \bS^2\}$ where $s$ is
the fixed nonzero constant in the definition of $Q_0$. From our
definitions in the next subsection, we shall see that the energy
well, $\Lambda_s$, corresponds to a set of uniaxial states.  Liquid
crystals satisfy the principle of frame indifference and are
macroscopically isotropic. As a consequence, $f_b$ is assumed to be
invariant with respect to orthogonal transformations, that is, we
require
\begin{equation}
f_b(RQR^t)=f_b(Q)\quad\text{for all }R\in O(3) \text{ and } Q \in \cS.
\end{equation}
Set
\begin{equation*}
\cS_0=\{Q\in \cS\colon Q_{i3}=Q_{3i}=0\quad\text{for }i=1,2\},
\end{equation*}

\begin{equation*}
\cA_0=\{Q(x,y)\in W^{1,2}(\Omega;\cS_0) {:}\quad Q=Q_0\text{ on
}\partial\Omega\},
\end{equation*}

and
\begin{equation*}
\cA=\{Q\in W^{1,2}(\Omega;\cS)\colon\quad Q=Q_0\text{ on
}\partial\Omega\}.
\end{equation*}

Our goal in this paper is to investigate minimizers for $F_\var$ in
$\cA_0$, and to analyze their behavior in the vanishing elastic
energy limit, $\var\to 0$. The relevance for doing this is that due
to the symmetries described above, these minimizers are critical
points (equilibria) for the energy $F_\var$ over the larger space
$\cA$, and thus satisfy the full set of Euler-Lagrange equations
with respect to variations in $\cA$.  (We prove this in Lemma 2.1.)
In addition, each $Q \in \cS$ is described in terms of an
orthonormal set of eigenvectors.  (See (1.8).)  For $Q \in \cS_0$,
we have
\begin{equation}
Q=s_1 \bm\otimes\bm+s_2\bm^\perp\otimes\bm^\perp-{1\over 3}\ (s_1+s_2)I
\end{equation}

for some real numbers $s_1$ and $s_2$, and $Q$ has an orthonormal basis of
eigenvectors of the form

\begin{eqnarray}
\{\bm,\bm^\perp,\be_3\}&\text{where } |\bm|=1,\ \bm=\langle
m_1,m_2,0\rangle,\\
&\text{ and } \bm^\perp=\langle -m_2,m_1,0\rangle,\nonumber
\end{eqnarray}
with eigenvalues

\begin{equation} \label{eigenvalues}
\lambda_1={1\over 3}(2s_1-s_2),\ \lambda_2={1\over
3}(2s_2-s_1),\ \lambda_3=-{1\over 3} (s_1+s_2).
\end{equation}

(See \cite{MN}.)  Thus the minimization problem of $F_\var$ over
$\cA_0$ models the behavior of a thin liquid crystal material
occupying $\Omega\times (-\eta,\eta)$ with its top and bottom
surfaces treated so as to fix $\be_3$ as a principal axis
(eigenvector of $Q$) of the liquid crystal molecules throughout the
body, with the other two principal axes (eigenvectors) in $\bR^2
\times \{0\}$, and boundary values on its side given by
$Q=Q_0(x,y)$.  The above problem includes a classic example from the
liquid crystal literature, in which
\begin{eqnarray}
f_b(Q)=f_b^0(Q)&=&\fa\ tr(Q^2)-{2\fb\over 3}\ tr(Q^3)+{\fc\over 2}
(tr (Q^2))^2+\fd\\
&=&\fa(\sum^3_{i=1}\lambda_i^2)-{2\fb\over 3}
(\sum^3_{i=1}\lambda_i^3)+{\fc\over
2}(\sum^3_{i=1}\lambda_i^2)^2+\fd \nonumber.
\end{eqnarray}
Indeed, taking $\fb,\fc>0$, $\fa<{\fb^2\over 27\fc}$, and  an
appropriate choice of $\fd$, we have $f_b^0\geq 0$ and $f_b^0(Q)=0$ if and
only
if $Q\in\Lambda_s$ where $s={1\over
4\fc}(\fb+\sqrt{\fb^2-24\fa\fc})$. (See \cite{MN}.)

\subsection{Definitions and Structural Assumptions}

 Our results require some structural assumptions on the bulk energy density $f_b$.
 In this section, we state these assumptions, along with some definitions and a
 change of variables in $\cA_0$ that will be needed to state our main results.

 It is well known (see \cite{MN}) that each $Q \in \cS$ has an orthonormal set
 of eigenvectors and can be written as
 \begin{equation}
Q=s_1 \bn\otimes\bn+s_2 \bk\otimes\bk-{1\over 3}\ (s_1+s_2)I
\end{equation}
where $\bn$ and $\bk$ are orthogonal unit vectors in $\mathbb R^3$; moreover,
the eigenvalues of $Q$ are given by the formula in (1.6).

 \begin{definition} Let $Q\in \cS$.  We say that $Q$ is
{\it isotropic} if all its eigenvalues are equal. (In this case, the structure
of $Q$ is that of a "normal" liquid.)

We say that $Q$ is {\it uniaxial} if exactly two of its eigenvalues
are equal.  (In this case, $Q$ has an axis of symmetry and its
structure is "rod-like" or "disk-like".)

We say that $Q$ is {\it biaxial} if all its eigenvalues are
distinct.  (In this case, there is no axis of complete rotational
symmetry for $Q$ and its structure is "board-like".)
\end{definition}

By formula (1.6) for the eigenvalues of $Q\in \cS$, it follows that $Q$ is
isotropic if and only if $s_1=s_2=0$ (and hence all eigenvalues are zero);
$Q$ is uniaxial if and only if one of the following three conditions hold:  $s_1=0$ and $s_2\neq 0$,  $s_2=0$ and $s_1\neq 0$, or $s_1=s_2\neq
0$ (and hence all eigenvalues are nonzero and exactly two of the eigenvalues
are equal).  Finally, $Q$ is biaxial for all other values of $s_1$ and $s_2$.

The above definition, when applied to a minimizer $Q_{\var}(\bx)$ of
$F_\var$ in $\cA$ or $\cA_0$, allows one to identify subregions of
$\Omega$ in which the liquid crystal material is in an isotropic,
uniaxial, or biaxial phase.  Note that $\Lambda_s\cap \cS_0$ is a
disconnected set of uniaxial states in $\cS_0$ with two connected
components: $\Lambda_s\cap \cS_0=\Lambda'_s\cup
\{s(\be_3\otimes\be_3-{1\over 3} I)\}$ where
$\Lambda'_s=\{s(\bm\otimes\bm-{1\over 3}I)\colon \bm=\langle
m_1,m_2,0\rangle,\ |\bm|=1\}$; also, the  boundary values $Q_0(x,y)$
are valued in $\Lambda'_s$.

\begin{definition}Let $\gamma\colon [0,1]\to\partial\Omega$ be a $C^3$
positively oriented parameterization of $\partial\Omega$ such that $\gamma$ is
one-to-one on $[0,1)$.  For $Q_0$ as assumed above, choose a unit vector field
$\tilde \bn_0(\bx)=\langle \tilde n_1(\bx),\tilde n_2(\bx),0\rangle$
defined on $\partial \Omega$ satisfying (1.1) such
that $\tilde \bn_0(\gamma(\cdot))\in C^1([0,1))$.  We
define the {\it degree} of $Q_0$ on $\partial\Omega$ by
\begin{equation*}
{1\over 2\pi}\ \int_0 ^1\tilde \bn_0
(\gamma(t))^\perp\cdot\frac{d\tilde\bn_0(\gamma(t))}{dt} \
dt\colon=\text{deg}~Q_0.
\end{equation*}
\end{definition}

Since $\lim\limits_{t\uparrow 1} \bn_0
(\gamma(t))=\pm\bn_0(\gamma(0))$ by (1.1) and the continuity of
$Q_0$, it follows  that deg~$ Q_0 ={k\over 2}$ for some $k\in \Bbb
Z$. Since we are interested in boundary conditions that correspond
to a liquid crystal with dislination-line type defects, we assume
that $k$ is nonzero, and thus without loss of generality, we shall
assume throught the paper that $k>0$. As $\var\downarrow 0$ the
effect of the bulk energy density $f_b$ becomes more pronounced and
minimizers tend to have their values located in a neighborhood of
$\Lambda_s\bigcap\cS_0$. Due to the boundary conditions, however,
this cannot happen throughout $\Omega$. We prove that the regions in
which minimizers, $Q_{\var}(x)$, of $F_{\var}$  take values outside
a neighborhood of $\Lambda'_s$ concentrate and quantize into $k$
small subdomains. For a subsequence as $\var_j\to 0$ these
subdomains tend to $k$ distinct points $\{a_1,\ldots,a_k\}$
representing the cross sections of the limiting disclination lines.

In \cite{SS} Schopohl and Sluckin carried out a numerical
investigation of equilibria for $F_\var$ in $\cA$. Their goal was to give
evidence that
equilibria are strongly biaxial near defects. They pointed out
that there is a subclass of equilibria which is contained in $\cA_0$,
and they developed simulations for these. This is the class of solutions
that we are studying here.

To state our main results, we will need the following linear change of variables
for the coefficients of each $Q
\in \cA_0$ in terms of unique functions $\bP=(p_1,p_2)$ and $r$:
\begin{equation}
Q=Q(\bP,r)=\begin{bmatrix}p_1+{r\over 2} &p_2 &0\\ p_2 &{r\over 2}-p_1 &0\\
0&0&-r\end{bmatrix}.
\end{equation}
From (1.1) and (1.4) each $Q\in\cA_0$ corresponds to a unique
$(\bP,r)\in W^{1,2}(\Omega;\bR^2)\times W^{1,2}(\Omega)$ satisfying
{$\bP|_{\partial\Omega}=\bP_0$, $r|_{\partial\Omega}=r_0$} where
$|\bP_0|={|s|\over 2}$, $r_0={s\over 3}$, and deg~$\bP_0=k=2 \text {
deg } Q_0$.  This
can be seen by writing (since $\bn_0\otimes\bn_0=(-\bn_0)\otimes(-\bn_0))$

\[
\bn_0(\gamma(t))=\pm\langle\cos\alpha(t),\sin\alpha(t),0\rangle
\]
for each $t$ in $[0,1)$, where $\langle
\cos\alpha(t),\sin\alpha(t),0\rangle=\tilde\bn_0(\gamma(t))$ and
$\alpha\in C^1([0,1)).$ Then using (1.1) and(1.4) we observe that

\[\bP_0(\gamma(t))={s\over2}\langle\cos {2\alpha(t)},\sin {2\alpha(t)}\rangle.
\]

 We may then recast our minimum problem by
considering the set
\begin{equation*}
A_0=\{(\bP,r)\in W^{1,2} (\Omega;\bR^2)\times W^{1,2}(\Omega)\colon
\bP=\bP_0 \text{ and }r={s\over 3}\text{   on }\partial\Omega\}.
\end{equation*}
The mapping $Q=Q(\bP,r): A_0 \to \cA_0$ is one-to-one and onto, and the
eigenvalues
for $Q(\bP,r)$ are $\lambda_1={r\over 2}+|\bP|$, $\lambda_2={r\over 2}-|\bP|$, $\lambda_3=-r$.
By (1.3) $f_b$ depends only on the invariants of
$Q$; since $trQ=0$, these are det $Q=({|\bP|^2-{r^2\over 4}})r$ and
$|Q|^2=2|\bP|^2+{3\over 2}\ r^2$. Thus $f_b(Q)=g_b(|\bP|^2,r)$ for
some function $g_b$. We prove in Section 2 that minimizing
$F_\var(Q)$ over $\cA_0$ is equivalent to minimizing
\begin{equation}
G_\var(\bP,r)=
\int_\Omega [g_e (\nabla\bP,\nabla r)+
\var^{-2} g_b (|\bP|^2,r)]\text{ for }(\bP,r)\in A_0,
\end{equation}
where $g_e (\nabla\bP,\nabla r)$ is defined by
\begin{eqnarray}
g_e&=&(L_1+{(L_2+L_3)\over 2})\ |\nabla\bP|^2+({3L_1\over 4}+{(L_2+L_3)\over 8})\ |\nabla r|^2\\
&+& {(L_2+L_3)\over 2}\ (p_{1x} r_x-p_{1y} r_y+r_x p_{2y}+r_y p_{2x})\nonumber\\
&+&|L_2+L_3| (p_{1x} p_{2y}-p_{1y} p_{2x}).\nonumber
\end{eqnarray}
This can be rewritten as
\begin{eqnarray}
g_e&=&L_1 (|\nabla\bP|^2+{3\over 4} |\nabla r|^2)\\
&+& {(L_2+L_3)\over 2} ((p_{1x}+{r_x\over
2}+p_{2y})^2+(p_{2x}-p_{1y}+{r_y\over 2})^2)\nonumber
\end{eqnarray}

\hskip4truein if $L_2+L_3\geq 0$,

\begin{eqnarray}
g_e&=&(L_1+L_2+L_3)(|\nabla\bP|^2+{3\over 4} |\nabla r|^2)\\
&-&{(L_2+L_3)\over 2}(({r_x\over 2}\!-\!p_{1x}-p_{2y})^2\!+\!(p_{2x}-p_{1y}-{r_y\over 2})^2\!+\!|\nabla r|^2)\nonumber
\end{eqnarray}

\hskip4truein if $0>L_2+L_3$.

The following structural conditions are assumed for $g_b(\fp,\fr)=g_b(|\bP|^2,r)$:

\begin{eqnarray}
\left\{\begin{array}{l}\text{ i)}\ g_b\in C^\infty([0,\infty) \times \mathbb R),
 g_b\geq 0 \text{ and } g_b({s^2\over4},{s\over 3})=0 {,}\\
\\
\text{  ii) }\text{ For some }m_1, m_2, m_3>0\\
~~~~|g_{b,\fp} (|\bP|^2,r) ||\bP|+|g_{b,\fr}(|\bP|^2,r)|\leq m_1(|\bP|^3+|r|^3)+m_2,\\
 ~~~~~~~~~m_3(|\bP|^4+|r|^4)-1\leq g_b (|\bP|^2,r) {,}\\
\\
\text{ iii) For some } \delta, m_4>0 \\
~~~~m_4((|\bP|^2-{s^2\over 4})^2+|r-{s\over 3}|^2)\leq g_b(|\bP|^2,r)\\
~~~~~~~~~~~\text{for }||\bP|-{|s|\over 2}|+|r-{s\over 3}|<\delta.
\end{array}\right.
\end{eqnarray}

Since $f_b(Q)=g_b(|\bP|^2,r)=g_b(\fp,\fr)$ under the change of
variables (1.9), these are additional assumptions on $f_b$.  From
(1.2), (1.12), and (1.13) we see that $g_e$ is a positive definite
quadratic. Thus $G_\var$ is strongly elliptic. It follows that
minimizers for $G_\var$ in $A_0$ exist and that the Euler--Lagrange
equation is a semi--linear elliptic system for which minimizers are
classical solutions $(C^\infty(\Omega)\bigcap C^2(\overline\Omega)$
in our case).  (See Theorem 2.2.)

In general the bulk energy well for $g_b(|\bP|^2,r)$ corresponds to $\{(\bP,r)\colon
g_b(|\bP|^2,r)=0\}$.  From our assumptions on $f_b$, the bulk energy well for
$f_b$ is
$\Lambda'_s\cup \{s(\be_3\otimes\be_3-{1\over 3} I)\}$. By the change of variables $Q \to
 (\bP,r)$,
$\Lambda'_s$ corresponds to  $\Gamma_s: =$ $\{(\bP,r)\colon
 |\bP|^2={s^2\over 4},r={s\over 3}\}$, and
 {$s(\be_3\otimes\be_3-{1\over 3} I)$} corresponds to
$(\bP,r)=(\mathbf 0,-{2s\over 3})$.  We note that the structural conditions
(1.14) only require that $\{g_b=0\}$ contains $\Gamma_s$ as in (i), that it is
bounded as in (ii), and that {$g_b$ has quadratic growth away from
$\Gamma_s$ as in (iii)}.

For the classic example, $f_b=f_b^0$ from the liquid crystal literature
(with coefficients $\fa,\fb,\fc$, and $\fd$ as described above (see (1.7)),
$f_b^0$ minimizes precisely on the uniaxial well $\Lambda_s$,

\begin{eqnarray*}
f_b^0(Q)&=&\fa(2|\bP|^2+{3\over 2}r^2)-2\fb r(|\bP|^2-{r^2\over 4})\\
&+&{\fc\over 2} (2|\bP|^2+{3\over 2} r^2)^2+\fd =\colon g_b^0
(|\bP|^2,r)
\end{eqnarray*}
for $Q\in {\cS}_0$ and $Q=Q(\bP,r)$, and one can easily show that
the structural assumptions (1.14) are satisfied for this example of
$g_b$.


\subsection{Main Results}

In this section we state our main results on the structure of
minimizers of the energy functional $F_{\var}(Q)$ over $\cA_0$,
using the fact that $Q$ is a minimizer of  $F_{\epsilon}$ in $\cA_0$
if and only if $(\bP,r)$ is a  minimizer of $G_{\epsilon}$ in $A_0$
and $Q=Q(\bP,r)$.

\medskip\noindent
{\it Theorem A. Let $\{(\bP_j,r_j)\}$ be a sequence of minimizers
for $\{G_{\var_j}\}$, respectively over $A_0$ such that
$\var_j\downarrow 0$. For ease of notation we consider $\bP_j$ as a
complex-valued function by identifying $\bR^2$ and $\bC$.  Then for
a subsequence $\{(\bP_{j'},r_{j'})\}$ there exists a harmonic
function $h\in C^2(\overline\Omega)$ and $k$ points
$\{a_1,\ldots,a_k\}\subset\Omega$ such that
\begin{eqnarray}
&&(|\bP_{j'}(\bx)|,r_{j'}(\bx))\to ({|s|\over 2},{s\over 3})\text{ in }C_{loc}
(\overline\Omega\backslash \{a_1,\ldots,a_k\}), \text{ and }\\
&&(\bP_{j'}(\bx),r_{j'}(\bx))\to (\bP^*(\bx),r^*(\bx))=({|s|\over 2}\
e^{i(h(\bx)+\sum^k_{\ell=1}\theta_\ell(\bx))}, {s\over 3})\nonumber
\end{eqnarray}
in $W_{loc}^{1,2}(\overline\Omega\backslash \{a_1,\ldots,a_k\}) \cap
C_{loc}(\overline\Omega\backslash \{a_1,\ldots,a_k\})$ and in
$C_{loc}^m(\Omega\backslash \{a_1,\ldots,a_k\})$ for all $m>0$,
where $\theta_\ell=\theta_\ell(\bx)$ denotes the polar angle of
$\bx$ with respect to the center $a_\ell$.  In particular, for each
sufficiently small $\rho > 0$, if $j'$ is sufficiently large,
setting $\Omega_\rho=\Omega\backslash\bigcup^k_{\ell=1} B_\rho
(a_\ell),$ we have
\begin{equation}
\bP_{j'}(\bx)=
|\bP_{j'}(\bx)|e^{i(h_{j'}(\bx)+\sum^k_{\ell=1}\theta_\ell(\bx))}\text{
in }\overline\Omega_\rho\,
\end{equation}

where $h_{j'}(\bx)$ is a function in $C^2(\overline\Omega_\rho)$ so
that $e^{ih_{j'}(\bx)}$ has degree zero on $\partial \Omega$, and
$\bP_{j'}$ has degree 1 about each of the k defects
$\{a_1,\ldots,a_k\}$.}

\bigskip

From the above result and the change of variables between $A_0$ and $\cA_0$, we
obtain:

\medskip\noindent
{\it Corollary A.  Let $\{Q_j\}$ be a sequence of minimizers of
$\{F_{\var_j}\}$, repectively over $\cA_0$ such that
$\var_j\downarrow 0$. Then for a subsequence of minimizers, we have
$Q _{j'}=Q(\bP_{j'}, r_{j'})$ where $\{\bP_{j'}, r_{j'}\} \subset
A_0$ satisfies Theorem A, and hence for each sufficiently small
$\rho > 0,$ if $j'$ is sufficiently large, we have:
\begin{eqnarray*}
Q_{j'}(\bx)&=& {s_{j'}}_1(\bx) (\bm_{j'}(\bx)\otimes
\bm_{j'}(\bx))+ {s_{j'}}_2(\bx) (\bm_{j'}^\perp(\bx)\otimes\bm_{j'}^\perp(\bx))\\
&-&{1\over 3} ({s_{j'}}_1(\bx)+{s_{j'}}_2(\bx)) I\qquad\text{in }\overline\Omega_\rho,
\end{eqnarray*}
where
\begin{eqnarray*}
\bm_{j'}(\bx)&=&\langle\cos({1\over 2} (h_{j'}(\bx)+\sum^k_{\ell=1}
\theta_\ell(\bx))),\sin({1\over 2}(h_{j'}(\bx)+\sum^k_{\ell=1}\theta_\ell(\bx))),0\rangle,\\
{s_{j'}}_1(\bx)&=& |\bP_{j'}(\bx)|+{3\over 2} r_{j'}(\bx),
\quad {s_{j'}}_2(\bx)={3\over
2}r_{j'}(\bx)-|\bP_{j'}(\bx)|,
\end{eqnarray*}

and $Q_{j'}$ has degree ${1\over 2}$ about each $a_\ell$.
(See Figure \ref{half}.)}

In particular, $Q_{j'}(\bx)$ converges to a uniaxial field
$Q^*(\bx)$ in $W_{loc}^{1,2}(\overline\Omega\backslash
\{a_1,\ldots,a_k\}) \cap C_{loc}(\overline\Omega\backslash
\{a_1,\ldots,a_k\})$ and in $C_{loc}^m(\Omega\backslash
\{a_1,\ldots,a_k\})$ for all $m>0$ as $j'\to\infty$, where

\[
Q^*(\bx) = s(\bm(\bx)\otimes\bm(\bx)-{1\over 3}\ I)\quad\text{in
}\Omega\backslash\{a_1,\ldots,a_k\}\text{ when } s > 0,
\]
and
\[
Q^*(\bx) =  s(\bm^\perp(\bx)\otimes \bm^\perp(\bx)-{1\over 3}\ I)\quad\text{in
}\Omega\backslash\{a_1,\ldots,a_k\}\text{ when } s < 0.
\]
Here
\begin{equation}
\bm(\bx)=\langle\cos ({1\over 2}\ (h(\bx)+\sum^k_{\ell=1}\theta_\ell(\bx))),
\sin({1\over 2}(h(\bx)+\sum^k_{i=1}\theta_\ell(\bx))),0\rangle
\end{equation}
for all $\bx$ in $\Omega\backslash\{a_1,\ldots,a_k\}$. Note that
$\bm_{j'}$ and $\bm$ are discontinuous while $Q_{j'}$ and $Q$ are continuous
on $\overline\Omega_\rho$.

The points $\{a_1,\ldots,a_k\}$ represent the cross sections of the
limiting disclination lines.  We prove in this paper that this set
of points minimizes a reduced energy $W(\bb)$ defined for
$\bb=(b_1,\ldots,b_k)\in\Omega^k$, which was introduced by Brezis,
Bethuel, and H\'elein in \cite{BBH} in connection with their
analysis of minimizing sequences $\{\bv_\var\}$ for the
Ginzburg--Landau energy
\begin{equation}
E_\var(\bv)={1\over 2}\int_\Omega [|\nabla\bv|^2+{1\over 2\var^2}
(1-|\bv|^2)^2]
\end{equation}
for $\bv\in\{\bw\in W^{1,2}(\Omega;\bR^2)\colon \bw=\bP_0/|\bP_0|
\text{ on } \partial\Omega\}.$  (The reduced energy $W(\bb)$ is
defined by equation (3.28).) More precisely, we have:

\begin{figure}[!htb]\label{half}
{\centerline{\includegraphics[angle=90, width=4in]{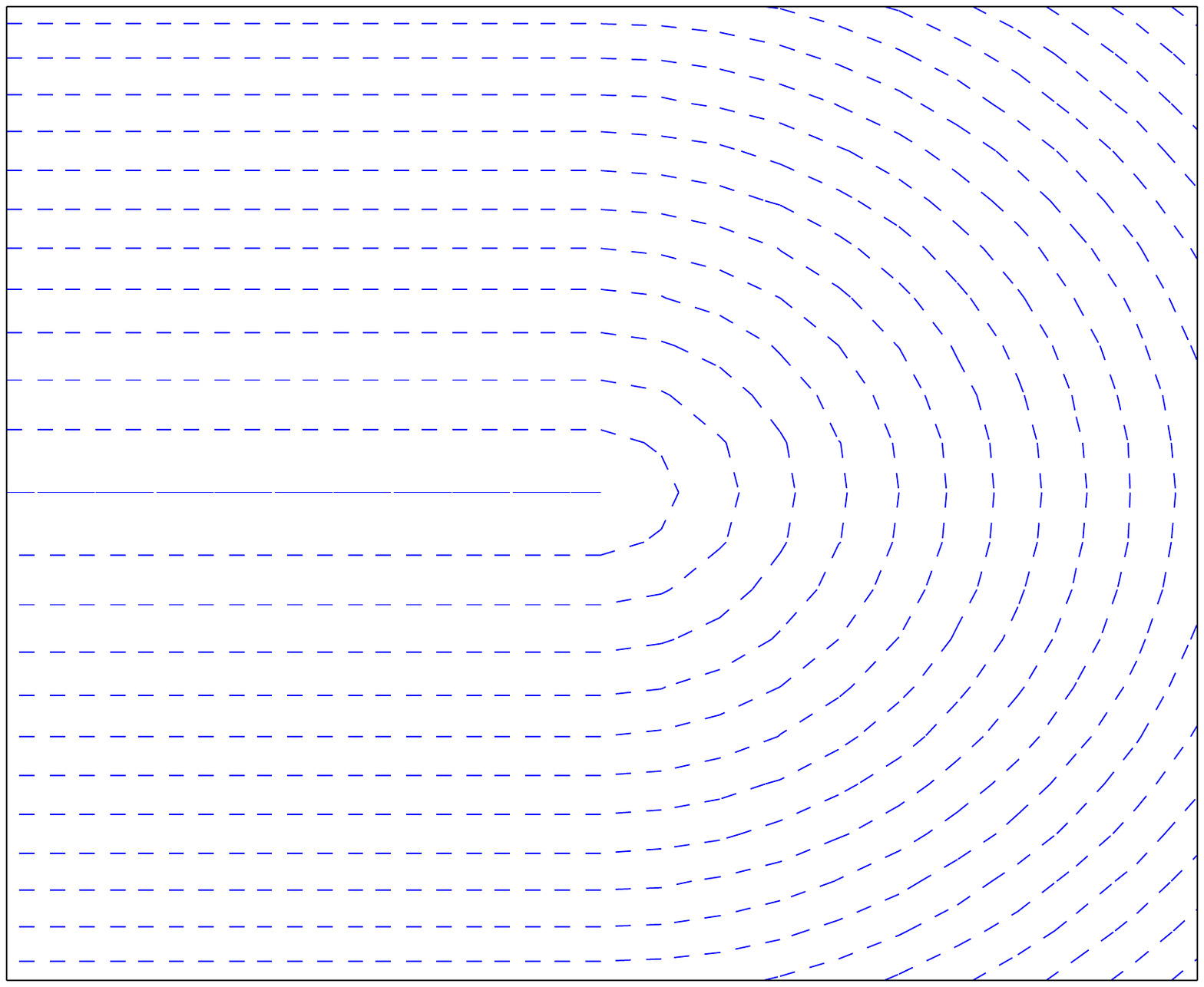}}}
\caption{${1\over 2}$ degree defect in a nematic texture.}
\end{figure}

\medskip\noindent
{\it Theorem B. Let $\{(\bP_j,r_j)\}$ be a sequence of minimizers
for $\{G_{\var_j}\}$ (or equivalently, let $\{Q_j\}$ be a sequence
of minimizers for $\{F_{\epsilon_j}\}$) for which $(a_1,\ldots,a_k)$
is a limiting configuration of defects as $\var_j\downarrow 0$ as
described in Theorem A. Then
\begin{equation*}
F_{\var_j}(Q_j)=G_{\var_j}(\bP_j,r_j)-(L_3-L_2+|L_3+L_2|){ s^2\pi
k\over4}.
\end{equation*}

Furthermore the reduced energy $W(\bb)$ for the limiting problem
minimizes at $\ba$ and we have
\begin{eqnarray*}
&\underset{j\to\infty}\lim&
[G_{\var_j}(\bP_j,r_j)-{{(2L_1+L_2+L_3)s^2\pi k}\over 4}\ln\
{1\over\var_j}]\\
&=&(2L_1+ L_2+L_3){{s^2}\over 4}W(\ba)+k\gamma.\\
\end{eqnarray*}
Here $\gamma$ is a fixed constant associated to the energy of each
defect core.}

The setting we study here gives a good description of
two-dimensional nematic behavior in flat films and thin layers.
Investigations from the physics literature of nematic textures in
flat and curved surfaces (thin shells) can be found in \cite{F},
\cite{LP}, \cite{N}, and \cite{VN}. In \cite{FS} Fatkullin and
Slastikov propose and investigate a model for two-dimensional
nematics (assuming that $L_1 > 0,$  $L_2=L_3=0$, and $Q=Q(\bP,r)$ is
in $\cA_0$ with $r(\bx) \equiv 0$) combining   Onsager- Maier-Saupe
and Landau de Gennes theories. This leads them to analyze a
variational problem closely related to the Ginzburg-Landau energy
(1.18).

Our last result describes how our work in this paper relates to earlier
investigations
of complex Ginzburg--Landau type functionals having
multiply-connected energy wells. The closest study in this respect
is \cite{HK} by Han and Kim in which they analyze the asymptotic
behavior for sequences of minimizers to the Chern-Simons-Higgs (CSH)
and the Maxwell-Chern-Simons-Higgs (MCSH) energies used to model
aspects of superconductivity.

For the (CSH) model one seeks (using our notation) minimizers
$\bP_\var$ to
\begin{equation}
C_\var(\bP)=\int_\Omega [{1\over
2}|\nabla\bP|^2+\var^{-2}|\bP|^2(1-|\bP|^2)^2]
\end{equation}
for $\bP\in B_0=\{\bv\in W^{1,2}(\Omega;\bR^2)\colon \bv=\bP_0$
 on $\partial\Omega\}$. Here $\bP_0\in C^3(\partial\Omega)$, $|\bP_0|=1$, and
 deg~$\bP_0=k>0$ with $k\in \Bbb N$.

For the (MCSH) model one seeks minimizers
$(\bP_{\var,q},r_{\var,q})$ to
\begin{equation}
C_{\var,q}(\bP,r)=\int_\Omega [{1\over 2}|\nabla\bP|^2+q^{-2}|\nabla
r|^2 +|\bP|^2r^2 +q^2(\var^{-1}(|\bP|^2-1)+r)^2]
\end{equation}
for $(\bP,r)\in B_0\times W^{1,2}_0(\Omega)$. The following two
results are from \cite{HK}:

i) For fixed $\var>0$, from any sequence of minimizers for (1.20)
with $q\to\infty$ one can find a subsequence
$\{(\bP_{\var,q_\ell},r_{\var,q_\ell})\}$  and a minimizer
$\bP_\var$ to (1.19) for which
\begin{eqnarray*}
\bP_{\var,q_\ell}\rightharpoonup\bP_\var \text{  and
}C_{\var,q_\ell}(\bP_{\var,q_\ell},r_{\var,q_\ell})\to
C_\var(\bP_\var)\text{ as } q_\ell\to\infty.
\end{eqnarray*}

ii) For fixed $q>0$, from any sequence of minimizers for (1.20) with
$\var \to 0$ there exists a subsequence
$\{(\bP_{\var_\ell,q},r_{\var_\ell,q})\}$, a point
$\ba^q=(a^q_1,\ldots,a^q_k)\in\Omega^k$, and a function $\bP_q^*$ as
in (1.15) so that $\bP_{\var_\ell,q}\to\bP_q^*$ in the sense of
Theorem A as $\var_\ell\to 0$ .

The functionals (1.10) and (1.20) are quite different. The bulk
energy well for $C_{\var,q}$ is
$\bS^1\times\{0\}\bigcup\{(\bf0,\var^{-1})\}$ and the second
component is eventually outside of any bounded set as $\var\to 0$.
This is in contrast to the bulk energy well for $G_\var$ which does
not vary with $\var$. The analysis in \cite{HK} is based on this
feature and it cannot be applied to (1.10).  Furthermore the bounds
in the estimates used to prove ii) diverge as $q\to\infty$ and so
they cannot be used to determine $\underset{\var\to 0}\lim
(\underset{q_\ell\to\infty}\lim
C_{\var,q_\ell}(\bP_{\var,q_\ell},r_{\var,q_\ell}))=\underset{\var\to
0}\lim C_\var(\bP_\var)$ or the limiting behavior of minimizers of
$C_\epsilon$, and this was left open. Our analysis however applies
to these issues directly. The same arguments we use to prove
Theorems A and B give the following result:

\medskip\noindent

{\it Theorem C. Let $\{\bP_\var\}$ be a sequence of minimizers for
(1.19) such that $\var\to 0$. Then there exists a subsequence
$\{\bP_{\var_\ell}\}$, a point $\ba=(a_1,\ldots,a_k)\in\Omega^k$,
and a function  $\bP^*$ as in (1.15) for which
$\bP_{\var_\ell}\to\bP^*$ in the sense of Theorem A. Moreover
$W(\cdot)$ minimizes at $\ba$ and
 \begin{eqnarray*}
\underset{\ell\to\infty}\lim[C_{\var_\ell}(\bP_{\var_\ell})-\pi
k\ln{1\over{\var_\ell}}] =W(\ba)+k\gamma
\end{eqnarray*} for a fixed constant
$\gamma$.}

\noindent Other related work is given in the papers \cite{KS1},
\cite{KS2}, and \cite{SY} in which the authors develop asymptotic
properties for the (CSH) energy using $\Gamma$--~convergence
techniques. This approach gives less detailed information than in
our setting. However, it is not restricted to sequences of
minimizers as {in} our case, and the authors apply it to more
general energies and scalings.

Our paper is organized as follows. In Section 2 we prove regularity
of minimizers and show that minimizers for $G_\var$ in $A_0$
correspond to a family of equilibria for $F_\var$ in $\cA$. In
Section 3 we prove Theorems A and B, developing the qualitative
features of minimizers for $G_\var$. Here we expand on
investigations of minimizers for the Ginzburg-Landau energy $E_\var$
(1.18) done by Brezis-Bethuel-H\'elein, Fanghua Lin, and Struwe.
(See \cite{BBH}, \cite{L1}, and \cite{St}). The energies $E_\var$
and $G_\var$ differ in two respects. The elastic term in the energy
density for $E_\var$ is the Dirichlet energy density, whereas for
$G_\var$ it is a coupled quadratic in $\nabla Q$. Secondly, the
energy well for the bulk energy density for $E_\var$ is $\bS^1$,
while the energy well for $G_\var$ is a bounded disconnected set
containing the ring $\Gamma_s$ as one of its components. The set
$\Gamma_s$ plays the same role as the energy well for $E_\var$. For
$\var$ small we prove that minimizers take their values near
$\Gamma_s$ except for an exceptional set contained in a neighborhood
of $k$ defects (vortices). In order to argue as has been done for
$E_\var$ we must first show that this exceptional set has small
measure. The results in Section 3 are proved assuming the a priori
estimate
\begin{equation}
\var^{-2}\int_\Omega g_b (|\bP_\var|^2,r_\var)\leq M
\end{equation}
for some constant $M<\infty$, for the family of equilibria
$\{(\bP_\var,r_\var)\colon0<\var<\var_1\}$  that are considered. In
Section 4 we prove, using a Pohozaev identity, that (1.21) is always
satisfied if $\Omega$ is a disk and $0<\var<1$. We then use this
result to establish (1.21) for the case in which $\Omega$ is a $C^3$
bounded simply connected domain and $\{(\bP_\var,r_\var)\}$ are
minimizers, where $\var_1$ depends on $s,L_1,L_2,L_3,\Omega,k,$ and
the constants in (1.14), and $M$ depends on these terms and in
addition on $\|\bP_0\|_{W^{1,2}(\partial\Omega)}$. Our approach for
this part is similar to one used by del Pino and Felmer in
\cite{dPF} in which they established the analogue of (1.21) for the
simpler energy (1.18).


\section{The Landau - de Gennes Energy}

By definition of $f_e$, we have
\begin{eqnarray*}
f_e(Q)&=&{L_1\over 2} |\nabla Q|^2+{(L_2+L_3)\over 2}|\text{div } Q|^2\\
&+& {L_3\over 2} (Q_{ij,k}Q_{ik,j}-Q_{ij,j}Q_{ik,k})
\end{eqnarray*}
where $\text{div }Q$ is the column vector whose $i$th entry is the divergence of the $i$th row of $Q$, $Q_{ij,j}$.
The last term in $f_e$ is a null--Lagrangian; its integral over $\Omega$ is constant on
\[
\cM=\{Q\in W^{1,2} (\Omega;\bR^{3\times 3})\colon Q=Q_0\text{ on }\partial\Omega\}
\]
and its first variation at any element of $\cM$ is zero. Set
$$f'_e(Q)={L_1\over 2} |\nabla Q|^2+{(L_2+L_3)\over 2}|\text{div}
Q|^2.$$
We say that $F_\var$ and $F'_\var=\int_\Omega
[f'_e+\var^{-2} f_b]$ are {\it equivalent} since their first
variations on $\cM$ agree, $\delta_VF_\var(Q)=\delta_VF_\var'(Q)$
where
\begin{eqnarray*}
&&\delta_V F (Q)=DF(Q) [V]:=\partial_t F (Q+tV)\text{ at
}t=0\\
&&\text{ for all }Q\in \cM \text{ and } V\in W^{1,2}_0 (\Omega;\bR^{3\times 3}).
\end{eqnarray*}

For $Q\in\cA$ we write
\begin{equation*}
Q(\bx)=\begin{bmatrix} z_1(\bx)&z_2(\bx)&z_4(\bx)\\ z_2(\bx)&z_3(\bx)&z_5(\bx)\\ z_4(\bx)&z_5(\bx)&-z_1(\bx)-z_3(\bx)\end{bmatrix}
\end{equation*}
and for $Q\in\cA_0$ we additionally have $z_4(\bx)=z_5(\bx)=0$. The
Euler--Lagrange equations for $F'_\var$ derived by variations in
$\cA(\cA_0)$ consist of the five (three) equations
$\delta_{z_\ell}F_\var'=0$ for $\ell=1,\ldots,5\ (\ell=1,2,3)$.

We can now show that an equilibrium with respect to variations in
$\cA_0$ is also an equilibrium with respect to variations in $\cA$.  We have:

\medskip\noindent
{\it Lemma 2.1.
 Let $Q\in \cA_0$  solve $\delta_{z_\ell}F_\var (Q)=0$ for $\ell=1,2,3$, then $\delta_{z_4} F_\var (Q)=\delta_{z_5} F_\var (Q)=0$ as well.}

\medskip\noindent
{\it Proof.} Since $f_b=\tilde f_b(\det Q, |Q|^2)$ it is easy to see that
$\partial_{z_4} f_b(\hat Q)=\partial_{z_5} f_b (\hat Q)=0$ for $\hat
Q\in \cS_0$. It follows directly that $\delta_{z_4} F'_\var
(Q)=\delta_{z_5} F'_\var (Q)=0$ for any $Q(\bx)\in\cA_0$.
\hfill$\square$

\medskip\noindent
{\it Theorem 2.2. For each $\var>0$ minimizers for $F_\var(Q)$ in
$\cA_0$ exist and are of class $C^\infty(\Omega)\bigcap
C^2(\overline\Omega)$.}

\medskip\noindent
{\it Proof.}
Recall that by (1.2), $L_1>0$ and $L_1+L_2+L_3>0$.  We consider two cases.

i)\ $L_2+L_3\geq 0$. From the discussion above we can work
with the energy $F'_\var$ instead of $F_\var$. Its energy density is the
sum of nonnegative
terms and $f'_e$ is a positive definite quadratic in $\nabla\bz$,
$\bz=(z_1,z_2,z_3)$. The first variation of $F'_\var$ in $\cA_0$
results in a semilinear elliptic system of three equations in three
unknowns. From standard elliptic theory (see \cite{G}) minimizers
for $F'_\var$ in $\cA_0$ exist, they are weak solutions to the
resulting elliptic system and they are classical
$(C^\infty(\Omega)\bigcap C^2(\overline\Omega))$.

ii)\ $0>L_2+L_3$. Let curl $Q$ denote the matrix-valued function
whose $i$th row is the curl of the $i$th row of $Q$. Then $|\nabla
Q|^2-|\text{div }Q|^2-|\text{curl
}Q|^2=(Q_{ij,k}Q_{ik,j}-Q_{ij,j}Q_{ik,k})$ is a null Lagrangian. As
a result, if we set
\[
f''_e(Q)={(L_1+L_2+L_3)\over 2}\ |\nabla Q|^2-{(L_2+L_3)\over 2}|\text{curl }Q|^2
\]
then $f_e-f''_e$ is a null Lagrangian, $f''_e$ is a positive
definite quadratic in $\nabla \bz$, and we can argue as in the
previous case.\hfill$\square$
\medskip

Setting $p_1=(z_1-z_3)/2, p_2=z_2$, and $r=z_1+z_3$ then $Q\in\cA_0$
is given in terms of $(\bP,r) \in A_0$ by (1.9). The minimum problem
for $F_\var$ in $\cA_0$ is recast as the minimum problem  for
$G_\var$ as defined in (1.10) in $A_0$, where $g_e$ as expressed in
(1.12) and (1.13) directly corresponds to $f'_e$ and $f''_e$ in
cases i) and ii) above, respectively. Moreover we have
$$g_\var(\bP,r)=f_\var(Q)-{(L_3-L_2+|L_3+L_2|)\over
4}(Q_{ij,k}Q_{ik,j}-Q_{ij,j}Q_{ik,k}).$$

\medskip\noindent
{\it Corollary 2.3. If $(\bP,r) \in A_0$ and $Q=Q(\bP,r)$ then
$$G_{\var}(\bP,r)=F_{\var}(Q)+(L_3-L_2+|L_3+L_2|){ s^2\pi
k\over 4}.$$}

\medskip\noindent
{\it Proof.} It suffices to evaluate
$\int_\Omega(Q_{ij,k}Q_{ik,j}-Q_{ij,j}Q_{ik,k}).$ As this is a
null-Lagrangian we are free to choose $(\bP,r) \in A_0$, and we set
$r={s\over3}$. It follows from (1.9) that
\begin{eqnarray*}
\int_\Omega(Q_{ij,k}Q_{ik,j}-Q_{ij,j}Q_{ik,k})&=&-4\int_\Omega(p_{1x}
p_{2y}-p_{1y} p_{2x})\\
&=&-4k|B_{s\over 2}(0)|=-s^2\pi k.
\end{eqnarray*}
\hfill$\square$

\medskip\noindent
{\it Corollary 2.4. Minimizers $(\bP_\var,r_\var)$ for $G_\var$ in
$A_0$ exist, they are of class $C^\infty(\Omega)\bigcap
C^2(\overline\Omega)$, and they correspond to minimizers for
$F_\var$ in $\cA_0$ by the relation (1.9).}


\section{The Asymptotic Problem}

By Theorem 2.2, equations (1.12)-(1.13), and our assumptions on $Q_0$,
it follows that minimizers $(\bP_\var,r_\var)$
for $G_\var$ in
$A_0$  are classical solutions to the boundary value problem
\begin{eqnarray}
&&\left\{\begin{array}{l}
\cL_1(\bP,r)\!:=\!-2 L_1\Delta p_1-(L_2+L_3)[\Delta p_1\!+\!{1\over 2}\! (r_{xx}\!-\!r_{yy})]\!=\!-{2p_1\over\var^2}\! g_{b,\fp}\\
\cL_2(\bP,r)\!:=\!-2L_1\Delta p_2-(L_2+L_3)[\Delta p_2+r_{xy}]\!=\!- {2p_2\over\var^2} g_{b,\fp}\\
\cL_3(\bP,r)\!:=\!-{3\over 2} L_1\Delta r-{(L_2\!+\!L_3)\over
2}[p_{1xx}\!-\!p_{1yy}\!+\!2 p_{2xy}\!+\!{1\over 2}\Delta
r]\!=\!-{1\over\var^2} g_{b,\fr}
\end{array}\right.\\
&&\hskip 10cm\text{in }\Omega,\nonumber\\
&&\qquad \qquad \text{ and }r={s\over 3},\qquad \bP=\bP_0\quad\text{ on }
\partial\Omega,
\end{eqnarray}
with $ |\bP_0|={|s|\over 2}$ on
$\partial \Omega$ and deg~$\bP_0=k>0.$\

Choose a finite covering $\cU$ of the $C^3$ manifold
$\overline\Omega$ by coordinate neighborhoods with uniformly bounded
$C^3$ structure, and a constant $\var_0$ in $(0,1)$ (depending only
on $\Omega$ and $\cU$) such that for all $x_0 \in\overline\Omega$,
$\overline B_{2\var_0}(x_0)$ is contained in a set in $\cU$.
Throughout this section we assume (1.21) holds for all minimizers
$\bz_\var=(\bP_\var,r_\var)$ for $G_\var$ in $A_0$ for all
$0<\var<\var_1\leq\var_0$, where $\var_1$ depends only on
$s,L_1,L_2,L_3,\Omega,k,$ and the constants in (1.14), and $M$
depends on these terms and in addition on
$\|\bP_0\|_{W^{1,2}(\partial\Omega)}$. This will be proved in
Section 4.

We begin this section by proving several a priori estimates, namely
Lemma 3.1 to Lemma 3.6, for solutions to (3.1) and (3.2) that
satisfy (1.21) for the above $M$ and $0<\var<\var_1$. These and
Proposition 3.7 to Corollary 3.13 will be applied to minimizers of
$G_\var$ to prove Theorems A and B at the end of this section.

In this section, unless otherwise stated, we denote by $C$ and $C_j$
positive constants depending at most on
$\bP_0,s,L_1,L_2,L_3,\Omega,$ and the constants in (1.14).
Additional dependence, e.g. on $M$, will be denoted by $C(M)$.

\medskip\noindent
{\it Lemma 3.1. Let $\bz_\var=(\bP_\var,r_\var)$ satisfy (1.21),
(3.1), and (3.2) for $0<\var<\var_1$. Then $|\bz_\var|$ and
$|\var\nabla\bz_\var|$ are uniformly bounded in $\overline\Omega$ by
a constant $C(M)$ independent of $\var$ for all $0<\var<\var_1.$}

\medskip\noindent
{\it Proof}. Let $\overline \bx\in\overline\Omega$ and let $\var \in
(0,\var_1)$. Set
\[
{\tilde \bz}(\by)=\bz_\var (\var \by+\overline \bx)\text{ for
}\by\in\overline{\tilde\Omega}=\{\by\colon \var \by+\overline
\bx\in\overline\Omega\}.
\]
Then in ${\tilde\Omega}$, ${\tilde \bz}$ satisfies the system obtained by
setting $\var=1$ in (3.1).
Let $\tilde B_r=B_r(0)\cap\overline{\tilde\Omega}$. From (1.21) and
the growth estimate (1.14) on $g_b$, we have
\begin{equation*}
\|\tilde{\bz}\|_{L^4(\tilde B_1)}\leq C(M) \qquad\text{for }
0<\var<\var_1.
\end{equation*}
Write (3.1) as $\cL\bz=\var^{-2}\bof(\bz)$, where $\bof(\bz)=[-2 p_1
g_{b,\fp},-2 p_2 g_{b,\fp},-g_{b,\fr}]^t$ and $\cL$ is {the second
order elliptic operator with constant coefficients}.  From
(1.21) and the $L^4$ estimate, we have
\[
\int_{\tilde B_1}|{\bof}({\tilde \bz})\cdot{\tilde \bz}|\leq C(M)
 \qquad\text{ for }0<\var<\var_1.
\]
In addition we have $\|{\tilde
\bz}|_{C^\ell(\partial\tilde\Omega)}\leq c_\ell$ for $0<\var<\var_1$
and  $\ell\leq 3$, where $c_\ell$ depends only on $\Omega$ and
$\bP_0$.

We use $\varphi^2(\tilde\bz-\psi)$ as a test function in (3.1) where
$\varphi$ is a cutoff function vanishing near $|\by|=1$, such that
$\varphi=1$ on $\tilde B_{3/4}$, and $\psi$ is a smooth function
equal to $\tilde\bz$ on $\partial\tilde\Omega$. The above
inequalities and elliptic estimates give $\|\tilde\bz\|_{1,2;\tilde
B_{3/4}}\leq C(M)$. This implies that $\mathbf f(\tilde\bz)\in
L^2(\tilde B_{3/4})$ and we see that $\|\tilde\bz\|_{2,2;\tilde
B_{5/8}}\leq C(M)$. Elliptic estimates imply that $\tilde\bz\in
W^{3,2}(\tilde B_{9/16})$ and by differentiating the equation we
obtain $\|\tilde\bz\|_{3,2;\tilde B_{9/16}}\leq C(M)$. It follows
that $\|\tilde\bz\|_{C^1(\overline{\tilde B}_{1/2})}\leq C(M)$
uniformly for $0<\var<\var_1$. The assertions then follow by scaling
back to $\bz_\var(\bx)$.\hfill$\square$

Set $\cO_\mu\colon=\{(\bP,r)\colon ||\bP|-{|s|\over 2}|+|r-{s\over
3}|\leq\mu\}$. Note that $\cO_0=\Gamma_s$. Below $\cH^n(E)$ denotes
the n-dimensional Hausdorff measure of $E$.

\medskip\noindent
{\it Lemma 3.2. Let $\bz_\var$ satisfy (1.21), (3.1), and (3.2). Set
$\cB(\var,\mu)=\{\bx\in\Omega\colon \bz_\var(\bx)\not\in\cO_\mu\}$,
$P_1(x,y)=x$, and $P_2(x,y)=y$. Let $0<\mu<\delta$ where $\delta$ is
given in (1.14). Then
\[
\cH^1(P_1 (\cB(\var,\mu))\leq C(\mu,M)\var \text{ and }
\cH^1(P_2(\cB(\var,\mu))\leq C(\mu,M)\var
\]
for all $0<\var<\var_1$.}

\medskip\noindent
{\it Proof}. Note that $\bz_\var(\bx)\in \cO_0$ for each
$\bx\in\partial\Omega$. Let $(x',y')\in\cB(\var,\mu)$, and set
$\ell_{x'}=\{(x',y)\colon y\in\bR\}$. Since this line intersects
$\partial\Omega$ there must exist $(x',y'')\in\ell_{x'}$ so that
$\bz(x',y'')\in\partial \cO_{\mu/2}$. It follows from Lemma 3.1 that
there is a $C_1(M)>0$ so that
\[
\bz(x',y)\in\cO_{3\mu/4}\backslash\cO_{\mu/4}\text{ for }
|y-y''|<C_1\var.
\]
From (1.14) then we see that there exists $C_2(\mu)>0$ so that
$C_2\var\leq\int_{\ell_{ x'}} g_b (|\bP_\var|^2,r_\var) d\cH^1(y)$.
Thus
\[
C_2\var\cH^1(P_1 (\cB(\var,\mu))\leq\int_\Omega
g_b(|\bP_\var|^2,r_\var)\leq\var^2M.
\]
The estimate for $P_2(\cB(\var,\mu))$ follows in the same manner.\hfill$\square$

Since $\cB(\var,\mu)\subset P_1(\cB(\var,\mu))\times
P_2(\cB(\var,\mu))$ for $\mu>0$ we have the following.

\medskip\noindent
{\it Corollary 3.3. Let $\bz_\var$ satisfy (1.21), (3.1), and (3.2).
For any $\mu\in(0,\delta)$  if $0<\var<\var_1$ then
$\cH^2(\cB(\var,\mu))\leq C(\mu,M)\var^2$.}

This estimate leads to a statement for all $\bx\in\Omega$. We  use
the fact that $(\bP_\var,r_\var)$ is bounded together with Corollary
3.3 for $\bx\in \cB(\var,\mu)$, and the growth estimate (1.14) for
$\bx\in\Omega\backslash\cB(\var,\mu)$ to get

\medskip\noindent
{\it Corollary 3.4. Let $\bz_\var$ satisfy (1.21), (3.1), and (3.2).
If \ $0<\var<\var_1$ then
\begin{equation}
\var^{-2}\int_\Omega ((r_\var (\bx)-{s\over 3})^2+ (|\bP_\var
(\bx)|^2-{s^2\over 4})^2)\leq C(M).
\end{equation}}

\medskip\noindent
{\it Lemma 3.5. Let $\bz_\var$ satisfy (1.21), (3.1), and (3.2). If
\ $0<\var<\var_1$ then
\[
\int_\Omega |\nabla r_\var|^2\leq C(M).
\]}

\medskip\noindent
{\it Proof}.
We first record an energy estimate for linear elliptic systems applied to (3.1) and (3.2),
\[
\|\bP_\var\|^2_{2,2;\Omega}+ \|r_\var\|^2_{2,2;\Omega}\leq
c_1(\var^{-4} (\|\bP_\var
g_{b,\fp}\|^2_{2;\Omega}+\|g_{b,\fr}\|^2_{2;\Omega})+\|\bP_0\|^2_{2,2;\partial\Omega})
\]
where $c_1$ depends on $L_1,L_2,L_3$ and $\Omega$. Since $g_b$
minimizes on $\cO_0$ we have
\begin{eqnarray}
&&|g_{b,\fp} (|\bP_\var(\bx)|^2, r_\var(\bx))|^2+|g_{b,\fr}(|\bP_\var(\bx)|^2, r_\var(\bx))|^2\\
&&\leq C((|\bP_\var (\bx)|^2-{s^2\over 4})^2+(r_\var (\bx)-{s\over
3})^2).\nonumber
\end{eqnarray}
Thus using (3.3) we find
\[
\|r_\var\|_{2,2;\Omega}^2\leq C(\var^{-2}+1).
\]
It then follows from this inequality and (3.3) that
\begin{eqnarray*}
\int_\Omega |\nabla r_\var|^2=-\int_\Omega (r_\var-{s\over 3})\Delta r_\var&\leq&\var^{-1} \|r_\var-{s\over 3}\|_{2;\Omega}\var\|r_\var\|_{2,2;\Omega}\\
&\leq& C(M).
\end{eqnarray*}
 \hfill{$\square$}

\medskip\noindent
{\it Lemma 3.6. There is a constant $\var_2\in (0,\var_1]$ depending
only on $\Omega$ and $k=\text{deg }\bP_0$, and a constant $C(M)$
independent of $\var$ so that if $(\bP_\var,r_\var)$ is a minimizer
for $G_\var$ in $A_0$ and $0<\var<\var_2$ then
\[
\int_\Omega |\nabla\bP_\var|^2\leq {s^2\over 4}\ 2\pi\ k\ln\
{1\over\var}+C(M).
\]}

\medskip\noindent
{\it Proof}. We first construct a comparison function for the energy
in (1.18). Choose  a set of distinct points
$\{b_1,\ldots,b_k\}\subset\Omega$, depending only on $\Omega$ and
$k$ such that
\begin{eqnarray*}
\min \{|b_n-b_\ell|,\text{ dist} (b_n,\partial\Omega); 1\leq
n,\ell\leq k,n\neq\ell\}=\overline\var
\end{eqnarray*}
is maximal. Define
\[
\bw_\var (\bx)=\prod^k_{\ell=1}\zeta({|\bx-b_\ell|\over\var } )\
{(\bx-b_\ell)\over |\bx-b_\ell| }\ e^{ij_\var(\bx)}
\]
where $\zeta(t)\in C^2(\bR)$ such that $\zeta(t)=0$ for $t\leq
{1\over 2}$, $\zeta(t)=1$ for $1\leq t$, and $j_\var(\cdot)$ is
harmonic in $\Omega$ such that $\bw_\var={\bP_0\over |\bP_0|}$ on
$\partial\Omega$ for $\var<\overline\var$. Then one has $E_\var
(\bw_\var)\leq\pi k\ln({1\over\var})+c_0$ for $0<\var<\overline\var$
where $E_\var$ is given in (1.18) and $c_0$ depends only on $\Omega$
and $\bP_0$. We next set $(\bw',r')=({|s|\over 2}\bw_\var,{s\over 3}
)\in A_0$ and use this as our comparison function for $G_\var$. Set
$\var_2=\text{ min }\{\overline\var,\var_1\}.$ Then for $\var\in
(0,\var_2]$, using (1.11) and (1.14) we find that
\begin{eqnarray*}
&&G_\var (\bw', r')\leq (L_1+ {L_2+L_3\over 2})\int_\Omega |\nabla \bw'|^2\\
&&\qquad\qquad\qquad + |L_1+L_3|\int_\Omega (w'_{1,x}
w'_{2,y}-w'_{1,y} w'_{2,x}) + C_1.
\end{eqnarray*}
The second integral on the right depends only on
$\bw'|_{\partial\Omega}$. Thus we get
\[
G_\var (\bw', r')\leq (L_1+ {L_2+L_3\over 2})\ {s^2\over 4}\ 2\pi k
\ln \ {1\over\var}+C_1.
\]
Next we use
\[
\int_\Omega g_e (\nabla\bP_\var, \nabla r_\var)\leq G_\var
(\bP_\var, r_\var)\leq G_\var (\bw', r').
\]
From (1.11) and suppressing the subscript $\var$ we see
\begin{eqnarray*}
&&(L_1+ {(L_2+L_3)\over 2})\int_\Omega |\nabla\bP|^2+{(L_2+L_3)\over 2}\int_\Omega (p_{1x} r_x-p_{1y}r_y+r_x p_{2y}+r_y p_{2x})\\
&+& |L_2+L_3| \int_\Omega (p_{1x} p_{2y} - p_{1y} p_{2x})\leq (2L_1+
L_2+L_3)\ {s^2\over 4}\pi k\ln({1\over\var})+C_1.
\end{eqnarray*}
Again the third integral is a constant depending on $\bP_0$. The
lemma will follow once we show that we can bound the second integral
appropriately. To do this we multiply the third equation in (3.1) by
$(r-{s\over3})$ and integrate over $\Omega$. We get using Lemma 3.5
that for $0<\var<\var_2:$
\begin{eqnarray*}
&&|{(L_2+L_3)\over 2}\int_\Omega (p_{1x} r_x-p_{1y} r_y +p_{2x} r_y+p_{2y} r_x)|\\
&\leq&\var^{-2}\int_\Omega |g_{b,\fr}|\cdot |r-{s\over 3}|+C_2(M)\\
&\leq& \var^{-2}\int_\Omega (|g_{b,\fr}|^2+ |r-{s\over
3}|^2)+C_2(M).
\end{eqnarray*}
Finally using (3.3) and (3.4) we see that {the last integral} is
bounded by a constant $C(M)$ independent of $\var$ for
$0<\var<\var_2$. \hfill $\square$

We are in a position to apply Lin's Structure Proposition, see
\cite{L3}. Significant parts of the proposition were also proved by
Sandier \cite{S} and Jerrard \cite{J}. Define

\begin{eqnarray*}
&&J_\var (\bv)=\int_\Omega j_\var(\bv),\text{ where}\\
&&j_\var(\bv)={1\over 2}[|\nabla\bv|^2 + {1\over 2\var^2}\
({s^2\over 4}-|\bv|^2)^2].
\end{eqnarray*}

\medskip\noindent
{{\it Proposition 3.7. For fixed $s\neq 0$ and a constant $K$
suppose that
\begin{eqnarray*}
&&\bP_\var\in \{\bv\in W^{1,2} (\Omega;\bR^2)\colon \bv=\bP_0\text{ on }\partial\Omega\}\text{ such that}\\
&&\bP_0\in W^{1,2}(\partial\Omega),~|\bP_0|={s\over 2},\text{ deg}~(\bP_0)=k>0,\\
&&J_\var (\bP_\var)\leq\pi\ {s^2\over 4}\ k\ln\ {1\over\var}+K
\end{eqnarray*}
where  $0<\var<\eta$. Fix $0<\alpha_0<min({1\over 8},{1\over
2(k+1)})$. There are positive constants $\eta_0\in(0,\eta)\text{ and
}\rho_0$ depending on $K,\Omega,\bP_0, \text{ and } \alpha_0$ so
that if $\var<\eta_0$ then for each $\bP_\var$ there are  points
$\{a_1^\var,\ldots,a_k^\var\}$
 for which
\begin{eqnarray*}
\min \{|a_n^\var-a_\ell^\var|,\text{ dist}
(a_n^\var,\partial\Omega); 1\leq n,\ell\leq k,n\neq\ell\}\geq\rho_0
\end{eqnarray*}
and constants $\alpha_m(\var)$, $\alpha_0\leq\alpha_m\leq 2\alpha_0$
for $1\leq m\leq k$ so that $|\bP_\var|\geq {|s|\over 2}$ on
$\partial B_m$ and deg~$(\bP_\var|_{\partial B_m})=1$ where
$B_m\colon =B_{\var^{\alpha_m}}(a_m^\var)$.

Moreover
\begin{equation}
\int_{\Omega\backslash\bigcup^k_{m=1} B_m} j_{\var}
(\bP_{\var})\leq\ {s^2\pi\over 4}\  {k\over k+1}\ \ln
{1\over\var}+c_1
\end{equation}

for some constant $ c_1=c_1(K,\Omega,\bP_0)$.

Furthermore for any sequence $\{\bP_{\var_\ell}\}$ with {$\var_\ell\downarrow 0$} there exists
a subsequence $\{\var_{\ell(q)}\}$, points $\{a_1,\ldots,a_k\}$, and a function $h(\bx)$ so that
\[
a_m^{\var_{\ell(q)}}\to a_m\text{ and }\bP_{\var_{\ell(q)}}\to\bP*=\prod_{m=1}^k\
{(\bx-a_m)\over |\bx-a_m|}\ e^{ih(\bx)}\ {|s|\over 2}
\]
as $q\to\infty$ where the convergence is strongly in $L^2(\Omega)$,
weakly in
\[
W_{loc}^{1,2}(\overline\Omega\backslash \{a_1,\ldots,a_k\}),\text{
and }\|h\|_{W^{1,2}(\Omega)}\leq c_2
\]
for some constant $c_2=c_2(K,\Omega,\bP_0)$.}}

 We take into account (1.21), Lemma 3.5,
Lemma 3.6 and apply the Proposition to a sequence of minimizers.

\medskip\noindent
{\it Lemma 3.8. Let $\{(\bP_\var,r_\var)\}$ be a sequence of
minimizers for $\{G_\var\}$ in $A_0$ such that $\var\downarrow 0$.
Then for a subsequence $\{(\bP_{\var_\ell},r_{\var_\ell})\}$ we have
$\bP_{\var_\ell}\to\bP^*$ as in Proposition 3.7 and
\[r_{\var_\ell}\rightharpoonup{s\over 3}\text{ in }W^{1,2}(\Omega).
\]}

The next two lemmas strengthen the notion of convergence using the
fact that we are working with a sequence of minimizers. Set
$\Omega_\rho=\Omega\backslash\bigcup\limits^k_{j=1}\ B_\rho (a_j)$.

\medskip\noindent
{\it Lemma 3.9. Let $\{(\bP_\ell,r_\ell)\}$ be a sequence of
minimizers for $\{G_{\var_\ell}\}$ in $A_0$ converging to
$(\bP^*,{s\over 3})$, in $L^2(\Omega)$ where $\bP^*(\bx)={|s|\over
2}\ \prod^k_{j=1}\ {(\bx-a_j)\over |\bx-a_j|}\ e^{ih(\bx)}$. Then
for each $0<\rho<{\rho_0\over2},$
\[
(\bP_\ell,r_\ell)\to(\bP^*,{s\over 3}) \text{ in }
W^{1,2}(\Omega_\rho)\ \text{ and }\ \underset{\var_\ell\to
0}\lim\var_\ell^{-2}\int_{\Omega_\rho}  g_b (|\bP_\ell|^2,r_\ell)=0.
\]
Moreover $\Delta h=0$ in  $\Omega$.}

\medskip\noindent
{\it Proof}. By Lemma 3.8 and an argument by contradiction we have
$(\bP_\ell,r_\ell)\rightharpoonup (\bP^*,{s\over 3})$ in $W^{1,2}
(\Omega_\rho)$ for each $\rho>0$ as above and  $(\bP_\ell,r_\ell)$
is a local minimizer for
\[
\int_{\Omega_\rho} [g_e (\nabla\bP,\nabla r)+\var_\ell^{-2} g_b (|\bP|^2,r)].
\]
To prove strong convergence it is enough to show that for each
$\overline \bx\in\overline\Omega\backslash \{a_1,\ldots,a_k\}$ there
exists a neighborhood $\cU_{\overline \bx}$ of $\overline \bx$ , on
which $(\bP_\ell,r_\ell)\to (\bP^*,{s\over 3})$ in
\newline
 $W^{1,2}
(\cU_{\overline \bx}\cap\Omega)$. We first consider the case
$\overline \bx\not\in\partial\Omega$ and take
$\overline{d}=\overline{d}(\overline \bx)>0$ such that $\overline
B_{2\overline{d}}=\overline B_{2\overline{d}}(\overline \bx)\subset
\Omega\backslash \{a_1,\ldots,a_k\}$. Then $\sum\limits^k_{j=1}
\theta_j(\bx)+h(\bx)$ is single valued here and we write
$\bP^*={|s|\over 2}\ e^{i\omega(\bx)}$ on $B_{2\overline d}$. From
Lemma 3.8 there exists {$C_0(M)<\infty$} independent of $\ell$ so
that
\[
\|\bP_\ell\|_{1,2;B_{2\overline d}}+\|r_\ell\|_{1,2;B_{2\overline
d}}\leq C_0(M).
\]
Thus for any subsequence $\{(\bP_{\ell_j},r_ {\ell_j})\}$ of
$\{(\bP_{\ell},r_ {\ell})\}$ (possibly after passing to a further
subsequence that we do not relabel) $d$ can be chosen, ${\overline
d}\leq d\leq 2{\overline d}$ so that
\begin{equation}
\int_{\partial B_{d}}[|\partial_\tau\bP_{\ell_j}|^2+|\partial_\tau
r_{\ell_j}|^2+\var_{\ell_j}^{-2}((|\bP_{\ell_j}|^2-{s^2\over
4})^2+(r_{\ell_j}-{s\over 3})^2)]\leq C_1(M)
\end{equation}
where  $\partial_{\tau}$  denotes the tangential derivative. Thus
$(|\bP_{\ell_j}|,r_{\ell_j})\to ({|s|\over 2},{s\over 3})$ uniformly
on $\partial B_d$ and $(\bP_{\ell_j},r_{\ell_j})\rightharpoonup
(\bP^*,{s\over 3})$ in $W^{1,2}(\partial B_d)$. Since
$\text{deg}~\bP^*|_{\partial B_d}=0$ it follows that
$\text{deg}~\bP_{\ell_j}|_{\partial B_d}=0$ for $j$ sufficiently
large and we can write $\bP_{\ell_j}(\bx)=|\bP_{\ell_j}(\bx)|
e^{i\omega_{\ell_j}(\bx)}$ for $\bx\in\partial B_{d}$. We define
$\tilde\omega_{\ell_j}(\bx)$ and $\tilde\omega(\bx)$ on $B_{d}$ as
the harmonic extensions of $\omega_{\ell_j}|_{\partial B_{d}}$ and
$\omega|_{\partial B_{d}}$ respectively. It follows that
\begin{equation}
\tilde\omega_{\ell_j}\rightharpoonup\omega\text{ in
}W^{1,2}(\partial B_{d})\text{ and
}\tilde\omega_{\ell_j}\to\tilde\omega\text{ in }W^{1,2}(B_d).
\end{equation}
The first limit follows from \cite{HKL} and the second follows from
elliptic regularity theory. We next construct comparison functions
\[
(\tilde\bP_{\ell_j},\tilde r_{\ell_j}) := (|\tilde\bP_{\ell_j}|
e^{i\tilde\omega_{\ell_j}},\tilde r_{\ell_j})\quad\text{ on }B_d
\]
such that $(\tilde\bP_{\ell_j},\tilde
r_{\ell_j})=(\bP_{\ell_j},r_{\ell_j})$ on $\partial B_{d}$.

This is done by setting
\[
(|\tilde\bP_{\ell_j}|,\tilde r_{\ell_j})=( {|s|\over 2},{s\over
3})\quad\text{ on }B_{d-\var_{\ell_j}},
\]
and for each $\theta$ define $(|\tilde\bP_{\ell_j}|,\tilde
r_{\ell_j})(|\bx|,\theta)$ to be linear for $d-\var_{\ell_j}\leq
|\bx|\leq d$. Then based on (3.6) and (3.7) it follows that
$(\tilde\bP_{\ell_j},\tilde r_{\ell_j})\to (\tilde\bP,\tilde
r)=({|s|\over 2}\ e^{i\tilde\omega},{s\over 3})$ in
$W^{1,2}(B_{d})$. Moreover
\[
\int_{B_d}g_e (\nabla\tilde\bP,0)=\lim_{j\to\infty}\
\int_{B_d}[g_e(\nabla\tilde\bP_{\ell_j},\nabla\tilde
r_{\ell_j})+\var_{\ell_j}^{-2} g_b (|\tilde\bP_{\ell_j}|^2,\tilde
r_{\ell_j})].
\]
From the minimality of $(\bP_{\ell_j},r_{\ell_j})$ and the weak
lower semicontinuity of $\int_{B_{d}}g_e$ we have
\begin{eqnarray}
\int_{B_d} g_e (\nabla\bP^*,0)&\leq&\limsup_{j\to\infty}\ \int_{B_d}
[g_e (\nabla\bP_{\ell_j},\nabla r_{\ell_j})+\var_{\ell_j}^{-2} g_b
(|\bP_{\ell_j}|^2, r_{\ell_j})]\\
&\leq&\int_{B_d} g_e(\nabla\tilde\bP,0)\nonumber.
\end{eqnarray}
From (1.11) it follows that $\int_{B_{d}} g_e(\nabla\bP,0)$
minimizes in the set $\{\bP={|s|\over 2}\ e^{if}\in W^{1,2}
(B_{d})\colon f=\omega$ on $\partial B_{d}\}$ if and only if $\Delta
f=0$ in $B_d$. Thus $\tilde\bP$ is the unique minimizer and
$\tilde\bP=\bP^*$ on $B_d$. From (1.12) and (1.13) we see that
$\int_{B_d} g_e (\nabla\bP,\nabla r)$ is the sum of weakly lower
semi--continuous integrals. We have shown that the sum is weakly
continuous on the sequence $\{(\bP_{\ell_j},r_{\ell_j})\}$. It
follows that each of its terms is weakly continuous on this sequence
as well. Thus $\int_{B_d}|\nabla\bP_{\ell_j}|^2\to\int_{B_d}
|\nabla\bP^*|^2$ and $\int_{B_d}|\nabla r_{\ell_j}|^2\to 0$ as
$j\to\infty$. Thus $(\bP_{\ell_j},r_{\ell_j})\to (\bP^*,{s\over 3})$
in $W^{1,2}(B_d)$ and as a result the full sequence
$(\bP_{\ell},r_{\ell})\to (\bP^*,{s\over 3})$ in
$W^{1,2}(B_{\overline d})$. A further consequence is that

\begin{equation*}
\lim_{\ell\to\infty}\ \var_\ell^{-2}\int_{B_{\overline d}}  g_b
(|\bP_\ell|^2, r_\ell)=0.
\end{equation*}

Moreover we have shown that {$\bP^*={|s|\over 2}\
e^{i(\sum_{j=1}^k\theta_j+h(\bx))}$} where $\Delta h=0$ in
$\Omega\backslash \{a_1,\ldots,a_k\}$. From Proposition 3.7 we have
that $h\in W^{1,2}(\Omega)$ and this implies that the singularities
are removable.

Lastly if $\overline \bx\in\partial\Omega$ we take a neighborhood
$\cU_{\overline \bx}$ and $d\in(0,\var_0)$ so that there exists a
smooth diffeomorphism defined on $B_{d}$ satisfying $\psi(\overline
\bx)=\overline \bx$ and
\[
\psi\colon B_{d}^+=\{ \by+\overline \bx\colon y_1^2+y_2^2 <d,\ y_2
\geq 0\}\underset{onto}\to\cU_{\overline \bx}.
\]
We can then carry out the radial construction of
$(|\tilde\bP_\ell|,\tilde r_\ell)$ in $B_{d}^+$, push this forward
to $\cU_{\overline \bx}$, and then argue as in the previous
case.\hfill$\square$

\medskip\noindent

 We next prove that $\{(|\bP_\ell|,r_\ell)\}$ {converges}
uniformly to $({|s|\over 2},{s\over 3})$ outside of a neighborhood
of $\{a_1,\ldots,a_k\}$. The proof is similar to that in {\cite{L1}
}Theorem A. This is possible since the density $g_e$ can be
expressed as the positive definite quadratic (1.12) or (1.13).

{\it Lemma 3.10. Let $(\bP_{\var_\ell},r_{\var_\ell})
=(\bP_\ell,r_\ell)$ be a {convergent} sequence of minimizers for
$\{G_{\var_\ell}\}$ in $A_0$ as in Lemma 3.9.  Given
$\rho\in(0,{\rho_0\over2})$ and $\mu\in(0,\delta)$
 there exists $\ell_0$ so that
\[
(\bP_\ell (\bx), r_\ell (\bx))\in\cO_\mu\text{ for all
}\bx\in\Omega_\rho\text{ and }\ell>\ell_0.
\]}

\medskip\noindent
{\it Proof}. Assume there exists $\bx_\ell\in\Omega_\rho$ such that
\[
(\bP_\ell (\bx_\ell), r_\ell (\bx_\ell))\not\in\cO_\mu\text{ for
}\ell\in \mathbb N.
\]

We consider two cases,
\begin{eqnarray*}
\text{i) dist}(\bx_\ell,\partial\Omega)\geq\var_\ell^{\alpha_0}\text{
for all }\ell,\\
\text{ii) dist}(\bx_\ell,\partial\Omega)<\var_\ell^{\alpha_0}\text{
for all }\ell,
\end{eqnarray*}
where we fix $0<\alpha_0<min({1\over 8},{1\over 2(k+2)})$  from Proposition 3.7.

We treat case $\text{ i)}$ first. Based on (3.3), Lemma 3.5, and
Lemma 3.6 we can select $\alpha_0<\alpha_\ell<2\alpha_0$ and set
$B'_\ell=B_{\var_\ell^{\alpha_\ell}}(\bx_\ell)$ so that
\begin{eqnarray}
\var_\ell^{\alpha_\ell}\int_{\partial B'_\ell}&&[(|\partial_\tau
\bP_\ell|^2+|\partial_\tau r_\ell|^2)+\var_\ell^{-2}
((|\bP_\ell|^2-{s^2\over 4})^2\\
&&+|r_\ell-{s\over 3}|^2)]\leq C_0\nonumber
\end{eqnarray}
for a fixed constant $C_0(M)$. Define
\begin{eqnarray*}
(\bP'_\ell (\by), r'_\ell (\by))=(\bP_\ell (\var_\ell^{\alpha_\ell}
\by+\bx_\ell), r_\ell(\var_\ell^{\alpha_\ell}\by+\bx_\ell))
\end{eqnarray*}
 for  $\by\in B_1(0)=B_1$. It follows that
\begin{equation}
|\nabla\bP'_\ell|+|\nabla r'_\ell|\leq C_1\var_\ell^{\alpha_\ell
-1} {,}
\end{equation}
for $C_1=C_1(M)$ and $(3.9)$ becomes
\begin{eqnarray}
\int_{\partial B_1}&&[(|\partial_\tau \bP'_\ell|^2+|\partial_\tau
r'_\ell|^2)+\var_\ell^{2(\alpha_\ell-1)}
((|\bP'_\ell|^2-{s^2\over 4})^2\\
&&+|r'_\ell-{s\over 3}|^2)]\leq C_0.\nonumber
\end{eqnarray}
We can assume that $B_\ell'\subset\Omega_{\rho/2}$ for each $\ell$.
Then using  \cite{L1} Lemma 1 for the first inequality and (3.5) for
the last we find
\begin{eqnarray*}
|\text{deg }\bP'_\ell|_{\partial B_1}|\ {s^2\pi \over 4}\ \
(1-\alpha_\ell)\ \ln\ {1\over\var_\ell}-C_2\ \leq\ \int_{B_1}
j_{{\var_\ell}^{(1-\alpha_\ell)}}
(\bP'_{\ell})\\
 =\int_{B_\ell'} j_{\var_\ell} (\bP_{\ell})\leq \ {s^2\pi
\over 4}\ {k\over k+1}\ \ln\ {1\over\var_\ell}+C_3.
\end{eqnarray*}

Since $\alpha_\ell<{1\over k+2}$ it follows that deg $
\bP'_\ell|_{\partial B_1}=0$ if $\ell$ is sufficiently large.
Moreover $(\bP'_\ell,r'_\ell)$ is a local minimizer for
\[
\int_{B_1}[g_e(\nabla\bP,\nabla r)+\var_\ell^{2(\alpha_\ell-1)} g_b
(|\bP|^2,r)].
\]
We can then construct comparison functions just as in Lemma 3.9, and
these lead as in the previous proof to
\begin{equation}
\lim\limits_{\ell\to\infty}\
\var_\ell^{2(\alpha_\ell-1)}\int_{B_1}g_b (|\bP'_\ell |^2, r'_\ell
)=0.
\end{equation}
On the other hand, using (3.11) it follows that
\[
(\bP'_\ell (\by), r'_\ell (\by))\in\cO_{\mu /4}\text{ for all
}\by\in\partial B_1\text{ and }\ell\geq\ell_1.
\]
 By hypothesis, we have  $(\bP'_\ell(0),r'_\ell(0))\not\in\cO_\mu$. Thus there
exists $\bz_{\ell}\in B_1$ with $(\bP'_\ell (\bz_\ell),
r'_\ell((\bz_\ell))\in\partial\cO_{3\mu/4}$. Using (1.14), (3.10),
and assuming $\mu<\delta$ (where $\delta$ is from (1.14)) we see
that there are {positive constants $C_4,\beta,$ depending in addition on $\mu$ and $M$,  so that
\[
g_b (|\bP'_\ell (\bx)|^2, r'_\ell(\bx))\geq\beta\text{ for }\bx\in
B_{C_4{\var_\ell}^{(1-\alpha_\ell)}}(\bz_\ell).
\]

Thus we conclude that
\[
\var_\ell^{2(\alpha_\ell-1)}\int_{B_1} g_b (|\bP'_\ell|^2,
r'_\ell)\geq C_5>0
\]
for a constant $C_5(\mu,M)>0$ and all $\ell$ sufficiently large.
This contradicts (3.12).

In case $\text{ii)}$ we consider $(\bP_{\var_\ell},r_{\var_\ell})$
on $B_{3\var_\ell^{\alpha_0}}(\by_\ell)\bigcap\Omega$ for
$\by_\ell\in\partial\Omega$ with
$|\bx_\ell-\by_\ell|\leq2\var_\ell^{\alpha_0}$. We can then flatten the
boundary to construct comparison functions as in the previous
lemma.\hfill$\square$

\medskip
In the next two lemmas we prove that if a sequence of minimizers
$\{(\bP_{\var_\ell},r_{\var_\ell})\}$ converges in
{$W_{loc}^{1,2}(\overline\Omega\backslash \{a_1,\ldots,a_k\})$} then
in fact it is bounded in {\newline$ W_{loc}^{j,2}(\Omega\backslash
\{a_1,\ldots,a_k\})$ }for all $j$. Our arguments are based on three
features, first that $\{(|\bP_{\var_\ell}|,r_{\var_\ell})\}$
converges uniformly to $({|s|\over 2},{s\over 3})$ on $K$ for each
$K\subset\subset\Omega\backslash \{a_1,\ldots,a_k\}$, second that
$({s^2\over 4},{s\over 3})$ is a nondegenerate minimum point for
$g_b$, and third that $g_e$ is strongly elliptic. A corresponding
result is proved for minimizing sequences to the Ginzburg--Landau
energy (1.18) in \cite{BBH}. In that case the Euler--Lagrange
equations are diagonal and the authors are able to apply estimates
for elliptic equations. Here our arguments rely only on $L^2$
estimates for elliptic systems.

\medskip\noindent
{\it Lemma 3.11. Let $\{(\bP_{\var_\ell},r_{\var_\ell})\}$ be a
sequence of minimizers for $\{G_{\var_\ell}\}$ in $A_0$ converging
in $W_{loc}^{1,2}(\overline\Omega\backslash \{a_1,\ldots,a_k\})$ as
$\var_\ell\to 0$. Then for $K\subset\subset\overline\Omega\backslash
\{a_1,\ldots,a_k\}$ there exist constants $\ell_0$ and $E$ so that
if $\ell\geq\ell_0$ then
\begin{equation*}
\|D^2(\bP_{\var_\ell},r_{\var_\ell})\|_{2;K}\leq E
\end{equation*}}

\medskip\noindent
{\it Proof}. It suffices to establish the estimate in a neighborhood
of each point in $\overline{\Omega}\backslash \{a_1,\ldots,a_k\}$.
We first consider the case of $x_0\in\Omega\backslash
\{a_1,\ldots,a_k\}$. Then
$\overline{B_{2d}}(\bx_0)\subset\Omega\backslash \{a_1,\ldots,a_k\}$
for some $d(\bx_0)\in (0,\var_0)$. Fixing $\bx_0$ and $\eta$,
$0<\eta<{|s|\over6}$, we take $d$ and $\ell_0$ so that
\begin{equation}
\int_{B_{2d}(\bx_0)} (|D\bP_{\var_\ell}|^2+|D r_{\var_\ell}|^2)<\eta
\end{equation}
and
\begin{equation}
\|\bP_{\var_\ell}|-{|s|\over 2}|+|r_{\var_\ell}-{s\over 3}|<
\eta\text{ on }B_{2d}(\bx_0)
\end{equation}
for all $\ell\geq\ell_0$.

 Let $\zeta\in C_c^2(B_{2d}(\bx_0))$ be
such that $\zeta=1$ on $B_d(\bx_0)$. We suppress the subscripts and
write $(\bP_{\var_\ell},r_{\var_\ell})=(\bP,r)$. Then multiplying
(3.1) by $-\partial_{x_j} (\zeta^2\partial_{x_j}(\bP,r))$, we get
using the strong ellipticity of the system that there exists a
constant $\Lambda(L_1,L_2,L_3)>0$ for which
\begin{eqnarray*}
&&\Lambda \|\zeta D\partial_{x_j}
(\bP,r)\|^2_{2;B_{2d}}+\var_\ell^{-2}\int_{B_{2d}} \zeta^2 [\cD^2
g_b]
(\partial_{x_j} (|\bP|^2,r))\cdot (\partial_{x_j} (|\bP|^2,r))\\
&&\leq C \||D\zeta|\partial_{x_j}
(\bP,r)\|^2_{2;B_{2d}}-\var_\ell^{-2}\int_{B_{2d}}2
g_{b,\fp}|\partial_{x_j}\bP|^2\zeta^2.
\end{eqnarray*}
Here $\cD g_b=(\partial_\fp g_b,\partial_\fr g_b)$ and $[\cD^2 g_b]$
is the Hessian of $g_b$. Using (1.14), (3.14) and taking $\eta$
sufficiently small we have
\[
\lambda\int_{B_{2d}}\zeta^2|\partial_{x_j}
(|\bP|^2,r)|^2\leq\int_{B_{2d}}\zeta^2 [\cD^2 g_b]\partial_{x_j}(
|\bP|^2,r)\cdot\partial_{x_j}(|\bP|^2,r)
\]
for some $\lambda>0$.

From  equations (3.1), using $|\bP|\geq {|s|\over 4}$ on $B_{2d}$,
we get
\[
\var^{-4}_\ell\int_{B_{2d}}\zeta^2
(g_{b,\fp}^2+g_{b,\fr}^2)=\var^{-4}_\ell\int_{B_{2d}} \zeta^2|\cD
g_p|^2\leq C\int_{B_{2d}}\zeta^2|D^2 (\bP,r)|^2.
\]
Thus we find
\begin{eqnarray}
&&\|\zeta D^2 (\bP,r)\|^2_{2;B_{2d}}+\var_\ell^{-4}\|\zeta\cD g_b\|^2_{2;B_{2d}}+\var_\ell^{-2}\|\zeta D(|\bP|^2,r)\|^2_{2;B_{2d}}\\
&\leq& C_0\int_{B_{2d}}\zeta^2 |D\bP|^4 + C_1\nonumber\\
&\leq& C_2\int_{B_{2d}}\zeta^2 |D^2\bP|^2\cdot\int_{B_{2d}}
|D\bP|^2+C_3.\nonumber
\end{eqnarray}
The last estimate follows by applying the Sobolev estimate
\begin{equation}
(\int_\Omega \varphi^2)^{1/2}\leq c\int_\Omega
(|D\varphi|+|\varphi|)
\end{equation}
with $\varphi=\zeta|D\bP|^2$ and $c=c(\Omega)$. Choosing $\eta$
small in (3.13) the first term on the right of (3.15) can be
absorbed into the left and the lemma is proved for the case of
$K\subset\subset\Omega\backslash\{a_1,\ldots,a_k\}$.

Assume next that $x_0\in\partial\Omega$ and $d<\var_0$, so that
$\overline{B_{2d}}(x_0)$ is contained in a coordinate patch in which
we can locally flatten $\partial\Omega$ near $x_0$. We consider the
special case where $\partial\Omega$ is already locally flat,
\begin{eqnarray*}
&&B_{2d}(x_0)\cap (\Omega\backslash\{a_1,\ldots,a_k\})=\\
&&B_{2d}^+(x_0)=\{(x_1,x_2)\colon (x_1-x_{01})^2+(x_2-x_{02})^2 <
4d^2 \text{ and } x_2 \geq x_{02}\}.
\end{eqnarray*}
Let $\zeta\in C_c^\infty (B_{2d} (x_0))$ such that $\zeta=1$ on
$B_d(x_0)$. Let $(\tilde\bP,\tilde r)\in W^{2,2}(\Omega)$ such that
$(\tilde\bP,\tilde r)=(\bP_0,{s\over3})$ on $\partial\Omega$. Again
suppressing subscripts, we multiply (3.1) by $\partial_{x_1}
(\zeta^2
\partial_{x_1} (\bP-\tilde\bP,r-\tilde r))$ and integrate by parts.
Then for any $0<\theta<1$ we get
\begin{eqnarray}
&&\Lambda\|\zeta D\partial_{x_1} (\bP,r)\|^2_{2;B_{2d}^+} \leq  C_1\| |D\zeta| |\partial_{x_1} (\bP,r)| \|^2_{2;B_{2d}^+}\\
&&+ \theta\var_\ell^{-4}\|\zeta\cD g_b\|^2_{2;B_{2d}^+}
+{1\over\theta}(\int_{B^+_{2d}} |\partial_{x_1} \bP|^4
\zeta^2+C_2).\nonumber
\end{eqnarray}
We next multiply (3.1) by
\[
-\partial_{x_2} (\zeta^2 \partial_{x_2} (\bP,r))=-\zeta^2
\partial_{x_2}^2 (\bP,r)-2\zeta \partial_{x_2} \zeta \partial_{x_2}
(\bP,r).
\]
Using the ellipticity of $\mathcal{L}$ we get
\begin{eqnarray}
&& {L_1\over 2}  \|\zeta^2\partial_{x_2}^2 (\bP,r)\|^2_{2;B_{2d}^+}-\Lambda_1 (\|\zeta^2 D\partial_{x_1} (\bP,r)\|^2_{2;B_{2d}^+}\\
&+& \| |D\zeta| |D (\bP,r)| \|^2_{2;B_{2d}^+})\nonumber\\
& \leq & -\int_{B_{2d}^+}\mathcal{L} (\bP,r)\cdot \partial_{x_2}
(\zeta^2 \partial_{x_2} (\bP,r))= I\nonumber
\end{eqnarray}
where $\Lambda_1=\Lambda_1 (L_1,L_2,L_3)$.

From (3.1) we have
\[
I=\int_{B_{2d}^+} [2p_1 g_{b,\fp}, 2p_2 g_{b,\fp},
g_{b,\fr}]^t\cdot\partial_{x_2} (\zeta^2\partial_{x_2} (\bP,r)).
\]
Here we integrate by parts. Since $g_b$ minimizes at
$(|\bP|^2,r)=(s^2,{s\over 3}),$ it follows that
$g_{b,\fp}=g_{b,\fr}=0$ on $\partial\Omega$. Thus the boundary term
will vanish and we find that
\begin{eqnarray}
I&=&-\var_\ell^{-2}\int_{B_{2d}^+} \partial_{x_2} [2 p_1 g_{b,\fp}, 2p_2 g_{b,\fp}, g_{b,\fr}]^t\zeta^2\partial_{x_2} (\bP,r)\\
&\leq& 2\var_\ell^{-2}\int_{B_{2d}^+} |g_{b,\fp}| |\partial_{x_2}
\bP|^2\zeta^2.\nonumber
\end{eqnarray}
Combining (3.17), (3.18), and (3.19) we see that there exists
$\Lambda_2(L_1,L_2,L_3)>0$ so that
\begin{eqnarray*}
&&\Lambda_2 (\|\zeta^2 D^2 (\bP,r)\|^2_{2;B_{2d}^+}+\var_\ell^{-4}\|\zeta\cD g_b\|^2_{2;B_{2d}^+}) \leq C_2\| |D\zeta| |D (\bP,r)| \|^2_{2;B_{2d}^+}\\
&&+\theta\var_\ell^{-4}\|\zeta\cD g_b\|^2_{2;B_{2d}^+}
+{1\over\theta}(\int_{B^+_{2d}} |D\bP|^4 \zeta^2+C_3).
\end{eqnarray*}
From this point the argument proceeds just as above. In the general
case one first flattens the boundary and analyzes the system in
local coordinates in the same manner. \hfill$\square$

\medskip\noindent
{\it Lemma 3.12. Let $\{(\bP_{\var_\ell},r_{\var_\ell})\}$ be the
sequence of minimizers for $\{G_{\var_\ell}\}$ from the previous
lemma. For each integer $j>2$ and set
$K\subset\subset\Omega\backslash \{a_1,\ldots,a_k\}$ there are
constants $E_j$ so that
\[
\|(\bP_{\var_\ell},r_{\var_\ell})\|_{j,2;K}\leq E_j\text{ for
}\ell\geq\ell_0.
\]}

\medskip\noindent
{\it Proof}. Choose $\eta<{|s|\over6}$ so that $[\cD^2g_b]\geq
\lambda I$ on $\cO_\eta$. We suppress the subscript $\var_\ell$ and
assume that $\ell\geq\ell_0$ where $\ell_0$ is from the previous
lemma. We further assume that  $d\in(0,\var_0)$ is sufficiently
small so that $\overline{B_{d}}(\bx_0)\subset\Omega\backslash
\{a_1,\ldots,a_k\}$ and so that (3.14) holds. Assume that there
exists a constant $E_q<\infty$ so that

\begin{eqnarray}
&&\|(\bP,r)\|^2_{q,2;B_{d}}+\var_\ell^{-2}\|(|\bP|^2-{{s^2}\over 4}, r-{s\over 3})\|^2_{q-1,2;B_{d}}\\
&+&\var_\ell^{-4}\| (g_{b,\fp},g_{b,r})\|^2_{q-2,2;B_{d}}\leq
E_q\nonumber
\end{eqnarray}
holds for $q=j-1$. We prove this estimate for $q=j$ where $E_{j-1}$
is replaced by a possibly larger constant, $E_j$ and $d$ by $d/2$.
Note that we already have (3.20) for $q=2$ from Lemma 3.11. Let
$\partial^\gamma$ be a derivative of order $j-1$ and $D^q$ be the
collection of all partial derivatives of order $q$. Let $\zeta\in
C_c^\infty (B_{d})$ be such that $\zeta=1$ on $B_{d/2}$. We use
$(-1)^{j-1}\partial^\gamma(\zeta^2\partial^\gamma(\bP,r))$ as a test
function in (3.1) and find
\begin{eqnarray}
&&\Lambda\|\zeta^2|D\partial^\gamma(\bP,r)| \|^2_{2;\Omega}\leq C\| |D\zeta|\partial^\gamma (\bP,r)\|^2_{2;\Omega}\\
&-&\var^{-2}_\ell\int_\Omega\zeta^2 \partial^\gamma(g_{b,\fp} 2\bP,
g_{b,r})\cdot\partial^\gamma(\bP,r)=I-\Pi\nonumber
\end{eqnarray}
From (3.20) we have $I\leq C_0(E_{j-1},d)$. We write
\begin{eqnarray}
&&\partial^\gamma (g_{b,\fp} 2\bP, g_{b,\fr})\cdot \partial^\gamma (\bP,r)=\partial^\gamma (g_{b,\fp},g_{b,\fr})\cdot (2\bP\cdot\partial^\gamma\bP,\partial^\gamma r)\\
&+&\sum_{|\alpha|\leq j-2\atop\alpha+\beta=\gamma} a_\alpha \partial^\alpha g_{b,\fp}\partial^\beta\bP\cdot\partial^\gamma\bP,\nonumber\\
&&2\bP\cdot\partial^\gamma\bP=\partial^\gamma
|\bP|^2+\sum_{\alpha+\beta=\gamma\atop 1\leq |\alpha|\leq j-2}
b_\alpha \partial^\alpha\bP\cdot\partial^\beta\bP,
\end{eqnarray}
and
\begin{eqnarray}
&&\partial^\gamma (g_{b,\fp},g_{b,\fr})=[\cD^2 g_b]\partial^\gamma (|\bP|^2,r)\\
&+&\sum_{\sum_{\alpha,\beta}(|\alpha|\ell_\alpha+|\beta|m_\beta)=j-1}c_{\alpha\beta}\prod_{|\alpha|\leq
j-2} (\partial^\alpha |\bP|^2)^{\ell_\alpha}\cdot \prod_{|\beta|\leq
j-2} (\partial^\beta r)^{m_\beta},\nonumber
\end{eqnarray}
where $a_\alpha,b_\alpha$ are constants, $\ell_0=m_0=0$, and
$c_{\alpha\beta}(\bx)=( c^1_{\alpha\beta} (\bx),
c^2_{\alpha\beta}(\bx))$ are bounded. Inserting (3.22), (3.23), and
(3.24) into the right side of (3.21) we have for $B_d=B_d(\bx_0)$:
\begin{eqnarray*}
II&=&\var_\ell^{-2}\int_{B_d}\zeta^2 [\cD^2 g_b]\partial^\gamma(|\bP|^2,r)\cdot\partial^\gamma (|\bP|^2,r) \\
&&+\var_\ell^{-2}\int_{B_d}\zeta^2\sum c_{\alpha\beta}\Pi(\partial^\alpha |\bP|^2)^{\ell_\alpha}(\Pi\partial^\beta r)^{m_\beta}\cdot\partial^\gamma(|\bP|^2,r) \nonumber\\
&&+\var_\ell^{-2}\int_{B_d}\zeta^2\partial^\gamma g_{b,\fp}(\sum b_\alpha\partial^\alpha\bP\cdot\partial^\beta\bP) \nonumber\\
&&+\var_\ell^{-2}\int_{B_d}\zeta^2 (\sum a_\alpha\partial^\alpha g_{b,\fp}\partial^\beta\bP\cdot\partial^\gamma\bP)\nonumber\\
&=&III+IV+V+VI.\nonumber
\end{eqnarray*}
Just as in Lemma 3.11 we have
\[
\lambda\var_\ell^{-2}\int_{B_d}\zeta^2 |\partial^\gamma
(|\bP|^2,r)|^2\leq III.
\]
From Sobolev's theorem the derivatives in IV of order less than
$j-2$ are bounded. It follows then for any $\theta>0$ that
\begin{eqnarray*}
|IV|&\leq& C_1\var_\ell^{-2}\int_{B_d}\zeta^2 (\sum_{t=1}^{j-2} |D^t(|\bP|^2,r)|^2) |\partial^\gamma (|\bP|^2,r)|\\
&\leq&\theta\var_\ell^{-4}\int_{B_d}\zeta^4 |D^{j-2}
(|\bP|^2,r)|^4+{C_2(E_{j-1},d)\over\theta}.
\end{eqnarray*}
Then using (3.16) and (3.20) we see
\[
|IV|\leq\theta
C_3(E_{j-1})\var_\ell^{-2}\int_{B_d}\zeta^2|D^{j-1}((|\bP|^2,r)|^2+{C_4(E_{j-1},d)\over\theta}.
\]
To estimate $|V|$ we write
$\partial^\gamma=\partial_{x'}\partial^{\gamma'}$ for some $x'$ and
integrate by parts to get
\[
|V|\leq\theta\var_\ell^{-4}\int_{B_d}\zeta^2| D^{j-2}
g_{b,\fp}|^2+\theta C_5(E_{j-1})\int_{B_d}\zeta^2| D^j \bP|^2+
{C_6(E_{j-1},d)\over\theta^2}.
\]
To bound $|VI|$ we first consider the terms with
$\alpha\ne\mathbf{0}$. For these $|\beta|<j-1$ and we see we can
bound these terms just as was done for $V$. The term with
$\alpha=\mathbf{0}$ can be bounded by
${C_7\over\theta}{g_{b,\fp}^2\over\var_\ell^{4}}+\theta |\zeta
D^{j-1}\bP|^4$. The integral of the first term over $B_d$ is bounded
from (3.20) and the second by $\theta C_8
(E_{j-1})\int_{B_d}\zeta^2| D^{j}\bP|^2+C_9(E_{j-1},d).$ Thus
\[
|VI|\leq C_{10}(E_{j-1})\theta \Big(\var_\ell^{-4}\int_{B_d}\zeta^2|
D^{j-2} g_{b,\fp}|^2+\int_{B_d}\zeta^2| D^{j}\bP|^2\Big)
+{C_{11}(E_{j-1},d)\over\theta}.
\]

\noindent Summing on $|\gamma|=j-1$ and collecting the estimates for
$III,\ldots,VI$ we find
\begin{eqnarray}
&&\Lambda\int_{B_d}\zeta^2|D^j(\bP,r)|^2+ \lambda(\var_\ell^{-2}\int_{B_d}\zeta^2| D^{j-1} (|\bP|^2,r)|^2\\
&\leq&\theta C_{12}(E_{j-1})\Big(\int_{B_d}\zeta^2|D^j (\bP,r)|^2+ \var_\ell^{-2}\int_{B_d} \zeta^2|D^{j-1} (|\bP|^2,r)|^2\nonumber\\
&+& \var_\ell^{-4} \int_{B_d}\zeta^2| D^{j-2} (g_{b,\fp},
g_{b,\fr})|^2\Big)+{C_{13}(E_{j-1},d)\over\theta^2}.\nonumber
\end{eqnarray}
From (3.1) we have $\var_\ell^{-2} (g_{b,\fp},g_{b,\fr})=-(
{\bP\over |\bP|^2} \cdot(\cL_1,\cL_2),\cL_3)(\bP,r)$. Using this,
the estimate $|\bP|\geq {|s|\over 4}$, and Sobolev's theorem we get
\[
\var_\ell^{-4}\int_{B_d}\zeta^2|D^{j-2} (g_{b,\fp},g_{b,\fr})|^2\leq
C_{14}(E_{j-1})\int_{B_d}\zeta^2|D^j (\bP,r)|^2+C_{15}(E_{j-1},d).
\]
Inserting this estimate into (3.25) and choosing $\theta$
sufficiently small we obtain (3.20) for $q=j$ and $d$ replaced by
$d/2$.\hfill$\square$

\medskip\noindent
{\it Corollary 3.13. Let $\{(\bP_{\var_\ell},r_{\var_\ell})\}$ be a
sequence of minimizers for $\{G_{\var_\ell}\}$ in $A_0$ converging
to $(\bP^*,r^*)$ in $W_{loc}^{1,2} (\overline{\Omega}\backslash
\{a_1,\ldots,a_k\})$. Then for each integer $m$
\[
(\bP_{\var_\ell},r_{\var_\ell})\to (\bP^*,r^*)\ {\rm in }\
C_{loc}(\overline\Omega\backslash \{a_1,\ldots,a_k\})
\]
and in $C_{loc}^m (\Omega\backslash \{a_1,\ldots,a_k\})$ as
$\ell\to\infty$.}

\medskip\noindent
{\it Proof of Theorem A}. Let $\{(\bP_\var,r_\var)\}$ be a sequence
of minimizers for $\{G_\var\}$ in $A_0$ for which (1.21) holds and
such that $\var\downarrow 0$. Then by applying Lemma 3.8 it follows
that there exists a subsequence
$\{(\bP_{\var_\ell},r_{\var_\ell})\}$ and points
$\{a_1,\ldots,a_k\}\subset\Omega$ so that
\begin{eqnarray*}
&&(\bP_{\var_\ell},r_{\var_\ell})\rightharpoonup ( {|s|\over 2}\ \prod^k_{j=1}\ {\bx-a_j\over |\bx-a_1|}\ e^{ih(\bx)}, {s\over 3})=(\bP^*,{s\over 3})\\
&&\text{in }W_{loc}^{1,2} (\overline\Omega\backslash
\{a_1,\ldots,a_k\})\times W^{1,2} (\Omega).
\end{eqnarray*}
By Lemma 3.10  for each $\rho\in (0,\var_0),$
$(|\bP_{\var_\ell}|,r_{\var_\ell})\to ({|s|\over 2},{s\over 3})$
uniformly on
$\overline\Omega_\rho=\overline\Omega\backslash\bigcup\limits^k_{j=1}
B_\rho (a_j)$, and from Lemma 3.9
\[
(\bP_{\var_\ell},r_{\var_\ell})\to (\bP^*,{s\over 3})\text{ in
}W^{1,2}(\Omega_\rho).
\]
Moreover $h(\bx)$ is harmonic in $\Omega$.

Finally, by applying Corollary 3.13 we see that
\[
(\bP_{\var_\ell},r_{\var_\ell})\to (\bP^*,{s\over 3})\text{ in
}C(\overline\Omega_\rho) \text{ and }C_{loc}^m (\Omega_\rho)
\]
for each integer $m$.\hfill$\square$

We need to establish several properties for the following minimum
problem in order to prove Theorem B. Let $\beta\in\bC$, $|\beta|=1$
and define
\begin{eqnarray}
&&L({\var\over\mu};\beta) := L({\var\over\mu},1;\beta)=L(\var,\mu;\beta)\nonumber\\
&=&\inf_{(\bv,r)\in\fA_\beta}\int_{B_\mu} [g_e (\nabla\bv,\nabla r)+\var^{-2} g_b (|\bv|^2,r)]\\
&+&(2L_1+L_2+L_3)\ {|s|\over 4}^2\pi \ln ({\var\over\mu})\nonumber
\end{eqnarray}
where
\[
\fA_\beta=\{(\bv,r)\in W^{1,2} (B_\mu)\colon \bv(\bx)={\beta
|s|\over 2}\ {x\over |x|} \text{ and } \ r(\bx)={s\over 3}\text{ for
}|x|=\mu\}.
\]

\medskip\noindent
{\it Lemma 3.14. $L(\tau;\beta)$ is independent of $\beta$ for all
$\beta\in\b C$ with $|\beta|=1$. Moreover $L(\tau):=L(\tau;\beta)$
is a nondecreasing function of $\tau$ for $\tau>0$ such that
$\gamma: =\lim\limits_{\tau\downarrow 0}\ L(\tau)
>-\infty$.}

\medskip\noindent
{\it Proof}. For any $T\in SO(2)$, consider the change of variables
by rotation, $\by=T\bx$ for $\bx\in B_1$ and set

\begin{equation*}
R=\begin{bmatrix}t_{11} &t_{12} &0\\ t_{21} &t_{22} &0\\
0&0&1\end{bmatrix}.
\end{equation*}

The energy density is frame indifferent and as such satisfies
\[
f_e(\nabla_\by\tilde Q(\by))+\tau^{-2}f_b(\tilde
Q(\by))=f_e(\nabla_\bx Q(\bx))+\tau^{-2}f_b( Q(\bx))
\]
where $\tilde Q(\by)=RQ(T^t\by)R^t.$ This translates into a
statement of invariance for $g_e$ and $g_b$,
\[
g_e(\nabla_\by\tilde \bP(\by),\nabla_\by\tilde
r(\by))+\tau^{-2}g_b(|\tilde \bP(\by)|^2, \tilde
r(\by))=g_e(\nabla_\bx \bP(\bx),\nabla_\bx r(\bx))+\tau^{-2}g_b(|
\bP(\bx)|^2, r(\bx))
\]
where $\tilde \bP(\by)=T^2\bP(T^t\by)$ and $\tilde
r(\by)=r(T^t\by).$ Let $\beta=\beta_1+i\beta_2.$ Then the boundary
condition for $\bP(\bx)$ as a vector in $\Bbb R^2$ reads as
$\bP_0(\bx)={|s|\over 2}Kx$ for $|\bx|=1$ where
\begin{equation*}
K=\begin{bmatrix}\beta_1 &-\beta_2\\\
\beta_2 &\beta_1\\
\end{bmatrix}.
\end{equation*}
Given $T\in SO(2)$ the boundary condition for $\tilde\bP(\by)$
becomes $\tilde\bP_0(\by)={|s|\over 2}T^2KT^t\by$ for $|\by|=1$. In
particular if we let $T=K^t$ we get $\tilde\bP_0(\by)={|s|\over
2}\by$ for $|\by|=1$. Thus the mapping
$(\bP,r)\in\fA_\beta\to(\tilde\bP,\tilde r)\in\fA_1$ is an isometry
such that $G_\tau(\bP,r)=G_\tau(\tilde\bP,\tilde r).$ In particular
we see that $L(\tau;\beta)=L(\tau;1)=L(\tau).$

 The monotonicity property of
$L(\tau)$ follows by the same argument for (1.18) given in [BBH],
Chapter 3. A lower bound $ \underline m$ for minimizers
 for the energy (1.18) with $\Omega=B_1$ is proved in [BBH], Chapter 5. Let $\bu_\var$
 be such a minimizer with $\bu_\var(\bx)={ \bx\over|\bx|}$ on $\partial B_1$. If
$(\bv_\tau,r_\tau)$ is a minimizer for (3.26) with $\mu=1$ and
$\var=\tau$ it follows that
\[
E_\tau({2\over|s|}\bv_\tau)\geq E_\tau( \bu_\tau)\geq -\pi\ln
(\tau)-\underline m.
\]
Thus using (3.3) we have
\[
{1\over 2}\int_{B_1} |\nabla\bv_\tau|^2\geq -{s^2\over 4}
\pi\ln(\tau)-\underline m'.
\]
The existence of a finite lower bound for $L(\tau)$ follows from
this and the estimates in the proof of Lemma 3.6.
 \hfill$\square$

\medskip\noindent
{\it Proof of  Theorem B}. The relation between $F_\var$ and
$G_\var$ is proved in Corollary 2.3. We establish the asymptotic
relation by arguing as in \cite{BBH}, Chapter 8. Let
\[
\Upsilon=\{\bb=(b,\ldots,b_k)\in\Omega^k\colon b_i\neq b_j\text{ if
}i\neq j\}
\]
and for $\bb\in \Upsilon$ set
\[
\bq_{b}(\bx)={|s|\over 2}\prod^k_{j=1}\ {(\bx-b_j)\over |\bx-b_j|}\
e^{i\bh_{\bb}(\bx)}
\]
where $\bh_{\bb}(\bx)$ is harmonic in $\Omega$ and is determined
(mod $2\pi$) by the condition $\bq_{\bb}=\bP_0$ on $\partial\Omega$.
From \cite{BBH}, Chapter 8 we have

\begin{eqnarray}
&&{1\over 2}\int_{\Omega\backslash\bigcup\limits^k_{j=1} B_\rho(b_j)} |\nabla\bq_{\bb}|^2=\\
&&{s^2\over 4}\ (\pi k\ln {1\over \rho}+W(\bb))+O(\rho)\text{ as
}\rho\to 0\nonumber
\end{eqnarray}
where $W(\bb)$ is the renormalized energy for (1.18) given in
\cite{BBH}. We express this using our notation. Set
$R(\bx)=\sum_{j=1}^{k} \ln |\bx-b_j|$ and
$\mathbf\tau=\mathbf\nu^\perp$ where $\mathbf \nu$ is the exterior
unit normal to $\partial\Omega$. Then

\begin{eqnarray}
W(\bb)&=&-\pi\sum_{\ell\ne j}\log |b_\ell-b_j|
+{1\over 2}\int_{\partial\Omega}R\partial_{\mathbf\nu} R\\
&&+\int_{\partial\Omega} h_{\bb}\partial_{\mathbf\tau} R + {1\over
2}\int_\Omega |\nabla h_{\bb}|^2.\nonumber
\end{eqnarray}

 Note that using (1.11) we have
\begin{eqnarray}
g_e(\nabla\bq_{\bb},0)&=&(L_1+{{L_2+L_3}\over 2})|\nabla\bq_{\bb}|^2\\
&+&|L_2+L_3| (q_{\bb1,x} q_{\bb2,y}-q_{\bb1,y}q_{\bb2,x})\nonumber
\end{eqnarray}
and that $q_{\bb1,x} q_{\bb2,y}-q_{\bb1,y}q_{\bb2,x}=0$ since
$|\bq_{\bb}|={|s|\over 2}$.

We next construct a comparison function for (1.10). Let $\bb\in
\Upsilon.$ Then for $0<\var_\ell<<\rho$ and for $\rho$ sufficiently
small (depending on $\Omega$ and $\bb$) we define
\begin{eqnarray*}
(\tilde\bP_{\var_\ell},\tilde r_{\var_\ell})=\left\{\begin{array}{l}
(\bq_{\bb}, s/3)\quad\text{ for
}\bx\in\Omega\backslash\bigcup^k_{j=1}\ B_\rho(b_j),\\
(\bv_0 (\bx-b_j), s/3)\quad\text{ for }\rho/2\leq |\bx-b_j|\leq\rho,
\\
(\bv_j ((\bx-b_j)),r_j (\bx-b_j))\quad\text{for }\bx\in B_{\rho/2}
(b_j).
\end{array} \right.
\end{eqnarray*}
Here $(\bv_j,r_j)$ minimizes
$\displaystyle\int_{B_{\rho/2}(0)}[g_e+\var_\ell^{-2} g_b]$ with
boundary conditions $({|s|\over 2}\ {\beta_j\bx\over |\bx|},{s\over
3})$ on $\partial B_{\rho/2}(0)$ and
$\beta_j=\prod\limits^k_{\ell=1\atop\ell\neq j}\ {(b_j-b_\ell)\over
|b_j-b_\ell|}\ e^{ih_{\bb}(b_j)}$.  The function $\bv_0$ is a
minimal harmonic map valued in $\{|\bv|=|{s\over 2}|\}$ such that
$\tilde\bP_{\var_\ell}$ is continuous. From Lemma 3.14 we have
\begin{eqnarray}
&&\int_{B_{\rho/2}(0)} [g_e (\nabla\bv_j,\nabla r_j)+\var_\ell^{-2} g_b (|\bv_j|^2,r_j)]\\
&=&(2L_1+ L_2+L_3){s^2\pi\over4}\ln ({\rho\over
2\var_\ell})+\gamma+o_{\var}(1)\nonumber
\end{eqnarray}
as $\var_\ell\to 0$.  Then from (3.27), (3.29), and Lemma 3.14 we
get
\begin{eqnarray*}
G_{\var_\ell} (\bP_{\var_\ell},r_{\var_\ell})&\leq& G(\tilde\bP_{\var_\ell},\tilde r_{\var_\ell})\\
&=&(2L_1+ L_2+L_3){{s^2}\over 4}(\pi k\ln ({1\over\var_\ell})+W(\bb))+k\gamma\\
&+&O(\rho)+o_\var (1).
\end{eqnarray*}
Let $\ba\in \Upsilon$ be a limiting configuration as in Theorem A.
Then from Lemma 3.9 and (3.27-30) we have
\begin{eqnarray*}
G_{\var_\ell}(\bP_{\var_\ell},r_{\var_\ell})&\geq&(2L_1+ L_2+L_3){{s^2}\over 4}(\pi k\ln ({1\over\var_\ell})+W(\ba))+k\gamma\\
&+&O(\rho)+o_\var (1).
\end{eqnarray*}
Just as in \cite{BBH}, choosing $\var_\ell=\var_\ell(\rho)<<\rho$
with $\rho\to 0$  we arrive at our assertion. It follows from these
two inequalities that $W$ minimizes at $\bb=\ba$ and that the limit
for $G_{\var_\ell} (\bP_{\var_\ell},r_{\var_\ell})$ as $\ell\to
\infty$ is established.  \hfill$\square$

\section{The Pohozaev Identity}

In this section we show that (1.21) always holds for minimizers of
$G_\var$ in $A_0$ if $\Omega$ is simply connected and
$0<\var<\var_1$ where $\var_1$ depends  on $s,L_1,L_2,L_3,\Omega,k,$
and the constants in (1.14), and $M$ depends on these terms and
$\|\bP_0\|_{W^{1,2}(\partial\Omega)}$ as well. We first prove (1.21)
for solutions to (3.1-2) in the case of a disk using the Pohozaev
identity.

\medskip\noindent
{\it Lemma 4.1. Let $(\bP,r)=(\bP_\var,r_\var)$ be a solution to
(3.1-2) where $\Omega=\Omega_R=B_R(0)$ and $0 <\var<1$. Then there
is a constant $C_0=C_0(R,L_1,L_2,L_3,\|\bP_0\|_{1,2;\partial
B_R},s)$ so that
\[
\var^{-2}\int_{B_R} g_b (|\bP|^2,r)\leq C_0.
\]}

\medskip\noindent
{\it Proof}. We multiply the system (3.1) by $-\nabla(p_1,p_2,r)\bx$
and integrate over $B_R$. We find
\begin{eqnarray}
&&0=\int_{B_R}\![(2L_1+L_2+L_3)\!(\Delta\bP\!\cdot\!\nabla\bP\!\cdot\!\bx)+({3L_1\over 2}\!+\!{L_2+L_3\over 4})\!\Delta r\nabla r\!\cdot\!\bx\\
&-&\var^{-2}\nabla g\cdot\bx]\nonumber\\
&+& {(L_2+L_3)\over 2}\int_{B_R} [2 r_{xy}\nabla p_2\cdot\bx+2{p_2}_{xy}\nabla r\cdot\bx]\nonumber\\
&+& {(L_2+L_3)\over 2}\int_{B_r} [(r_{xx}-r_{yy})\nabla p_1\cdot\bx+ ({p_1}_{xx}-{p_1}_{yy})\nabla r\cdot\bx]\nonumber\\
&=:& I+{(L_2+L_3)\over 2}\ II+{(L_2+L_3)\over 2}\ III.\nonumber
\end{eqnarray}
We can calculate $I$ as in \cite{BBH}, Chapter 3,
\begin{eqnarray}
I&=&R (L_1+ {(L_2+L_3)\over 2})\int_{\partial B_R} (|\bP_\nu|^2-|\bP_\tau|^2)\\
&+& R({3L_1\over 4}+{(L_2+L_3)\over 8} )\int_{\partial B_R} (|
r_\nu|^2-|r_\tau|^2)+2\var^{-2}\int_{B_R} g_b.\nonumber
\end{eqnarray}
Here $\bP_\tau$ and $r_\tau$ are tangential derivatives. Note that
$r_\tau=0$ and $\bP_\tau=\bP_{0\tau}$ on $\partial B_R$. To
calculate II we write
\begin{eqnarray*}
&&\int_{B_R} r_{xy}\nabla p_2\cdot\bx=\int_{B_R} (r_{xy}x {p_2}_x+ r_{xy} y {p_2}_y)\\
&=& -\int_{B_R} (x r_x {p_2}_{xy} + y r_y {p_2}_{xy})\\
&+& {1\over R}\int_{\partial B_R} {xy} ({p_2}_x r_x + {p_2}_y r_y).
\end{eqnarray*}
Using this and the fact that $r_\tau=0$ on $\partial B_R$ we get
\[
II={2\over R}\ \int_{\partial B_R} xy{p_2}_\nu r_\nu.
\]
To calculate III we change variables, $x'=(x-y)/\sqrt{2}$,
$y'=(x+y)/\sqrt{2}$. Then
\[
III=2\int_{B_R} (r_{x' y'}\nabla p_1\cdot\bx+ {p_1}_{x' y'} \nabla
r\cdot\bx)={2\over R}\ \int_{\partial B_R} x' y' {p_1}_\nu r_\nu.
\]
Writing $(x,y)=(R\cos\theta,R\sin\theta)$ then it follows that
$(x',y')=(R\cos(\theta+{\pi\over 4}), R\sin (\theta+{\pi\over 4}))$.
Thus $II+III=R\int_{\partial B_R} r_\nu (\cos 2\theta,\sin
2\theta)\cdot \bP_\nu$. Finally we see that
\begin{equation}
\bigg|\left({L_2+L_3\over 2}\right) (II+III)\bigg|\leq R\
{|L_2+L_3|\over 2}\ \left(\int_{\partial B_R}( {|r_\nu|\over
4}^2+|\bP_\nu|^2)\right).
\end{equation}
Thus using (4.1), (4.2) and (4.3) with (1.2) we get
\begin{eqnarray*}
&& R\left(L_1+{(L_2+L_3)\over 2}\right)\int_{\partial B_R} |{\bP_0}_\tau|^2\\
&\geq& R(L_1+ {L_2+L_3\over 2}-{|L_2+L_3|\over 2})\int_{\partial B_R} |\bP_\nu|^2\\
&+&R \left( {3L_1\over 4}+{(L_2+L_3)\over 8}-{ |L_2+L_3|\over 8}\right)\int_{\partial B_R} |r_\nu|^2\\
&+&2\var^{-2}\int_{B_R} g_b\geq 2\var^{-2}\int_{B_R} g_b.\\
\end{eqnarray*}
\rightline{$\square$}

\medskip\noindent
{\it Lemma 4.2. Let $\Omega$ be a $C^3$ bounded simply connected
domain in $\bR^2$. There is a constant $0<\var_1\leq\var_0$ such
that if $(\bP,r)=(\bP_\var,r_\var)$ is a minimizer for $G_\var$ in
$A_0$ and $0<\var<\var_1$, then
\[
\var^{-2}\int_\Omega g_b (|\bP|^2,r)\leq M.
\]
Here  $\var_1$ depends on  $s,L_1,L_2,L_3,\Omega,k,$ and the
constants in (1.14) and $M$ depends on these terms and
$\|\bP_0\|_{W^{1,2}(\partial\Omega)}$.}
\medskip\noindent

{\it Proof}. Set $R=2 (\text{diam}(\Omega))$ and assume that
$0\in\Omega$. We construct an extension of $\bP.$ Let $\hat\bP\in
W^{1,2}( B_R(0)\backslash \Omega)$ valued in $\{|\hat\bP|={|s|\over
2}\}$ and such that $\hat\bP$ is a minimal harmonic map satisfying
$\hat\bP=\bP_0$ on $\partial\Omega$ and
$\hat\bP(\bx)={|s|\over2}({\bx\over|\bx|})^k$ on $\partial B_R(0)$.
Note that $\|\hat\bP\|_{1,2;B_R(0)\backslash \Omega}\leq
C\|\bP_0\|_{1,2;\partial\Omega}.$ Set

  \begin{eqnarray*}
(\bP',r')=\left\{\begin{array}{l} (\bP,r) \text{ for }\bx\in\Omega,\\
(\hat{\bP},{s\over 3}) \text{ for }\bx\in B_R\setminus\Omega.
\end{array} \right.
\end{eqnarray*}

Let $\tilde G_\var=\displaystyle\int_{B_R} [g_e+{1\over 2\var^2}\
g_b]$, and let $(\tilde\bP,\tilde r)$ be a minimizer for $\tilde
G_\var$ such that $(\tilde\bP,\tilde r)=(\hat\bP,{s\over 3})$ on
$\partial B_R$. We can apply Lemma 4.1 (with $\var$ replaced by
$\sqrt{2}\var$) and the results from Section 3 to $\tilde G_\var$
and $(\tilde\bP,\tilde r)$ for the case of $\Omega=B_R$. In
particular { from the proof of Theorem B} there are  constants $C_1$
and $0<\eta_1<1$, depending on $s,L_1,L_2,L_3,\Omega,k,$ and the
constants in (1.14) so that
\[
(2L_1+ L_2+L_3){s^2\over 4}\ \pi k\ln {1\over\var}-C_1 \leq\tilde
G_\var (\tilde\bP,\tilde r)\leq \tilde G_\var (\bP',r')
\]
for all $0<\var<\eta_1$. Note that
\begin{eqnarray*}
\tilde G_\var (\bP',r')&=&\int_\Omega [g_e (\nabla\bP,\nabla r)+{1\over 2\var^2}\ g_b (|\bP|^2,r)]\\
&+&\int_{B_R\backslash\Omega} g_e (\nabla\hat\bP,0)\\
&=& G_\var(\bP,r)-{1\over 2\var^2}\ \int_\Omega g_b(|\bP|^2,r)+C_2
\end{eqnarray*}
where $C_2$ depends only on $\|\bP_0\|_{1,2;\partial\Omega}$ and the
constants in (1.14). Thus
\begin{equation}
(2L_1+ L_2+L_3){s^2\over 4}\pi k\ln {1\over\var}+{1\over
2\var^2}\int_\Omega g_b(|\bP|^2,r)\leq G_\var (\bP,r)+C_1+C_2.
\end{equation}

Next we consider the comparison map $(\bw',r')$ constructed in Lemma
3.6 defined for $\var<\overline\var=\eta_2$. Since $(\bP,r)$ is a
minimizer for $G_\var$ we get
\[
G_\var(\bP,r)\leq G_\var(\bw',r')\leq (2L_1+ L_2+L_3){s^2\over 4}\pi
k\ln {1\over\var}+C_3
\]
for all $\var<\overline\var=\eta_1$ where $C_3$ depends only on
$\bP_0,\Omega,L_1,L_2,L_3,$ and the constants in (1.14). It follows
from this and (4.4) that
\[
\var^{-2}\int_\Omega g_b (|\bP|^2,r)\leq 2( C_1+C_2+C_3)=\colon M.
\]
for all $0<\var<\var_1=\text{min}\{\eta_1,\eta_2,\var_0\}$.

\rightline{$\square$}

\end{document}